\documentclass[11pt]{amsart}
\usepackage{epsfig}
\usepackage{psfrag}
\newtheorem{theorem}{Theorem}
\newtheorem{proposition}{Proposition}
\newtheorem{lemma}{Lemma}
\newtheorem{corollary}{Corollary}
\newtheorem{definition}{Definition}
\numberwithin{equation}{section}
 \numberwithin{theorem}{section}
 \numberwithin{lemma}{section}
\numberwithin{proposition}{section}
 \numberwithin{corollary}{section}
\numberwithin{definition}{section}
\def\R{\bf R}
\def\e{\epsilon}

\def\al{\aligned}
\def\eal{\endaligned}
\def\M{{\bf M}}
\def\be{\begin{equation}}
\def\ee{\end{equation}}
\def\lab{\label}

\begin{document}

\tracingpages 1
\title[extremal]{\bf Extremal of Log Sobolev inequality and $W$ entropy on
noncompact manifolds}
\author{Qi S. Zhang}
\address{ Department of
Mathematics, University of California, Riverside, CA 92521, USA }
\date{May 2011}

\begin{abstract}

Let $\M$ be a complete, connected noncompact manifold with bounded
geometry. Under a condition near infinity, we prove that the Log
Sobolev functional (\ref{logfanhan})  has an extremal function
decaying exponentially near infinity. We also prove that an
extremal function may not exist if the condition is violated. This
result has the following consequences. 1. It seems to give the
first example of connected, complete manifolds with bounded
geometry where a standard Log Sobolev inequality does not have an
extremal.
 2. It gives a negative answer to the open question on the existence of extremal of Perelman's
 $W$ entropy in the noncompact case, which was stipulated by Perelman \cite{P:1} p9, 3.2 Remark. 3.
It  helps to prove, in some cases,  that noncompact shrinking breathers of Ricci flow are gradient shrinking solitons.

\end{abstract}
\maketitle \tableofcontents
\section{Introduction}

The main purpose of the paper is to give a counter example to the
old question on existence of extremals of a standard Log Sobolev
inequality (or its recent reincarnation in the form of Perelman's
$W$ entropy) on noncompact manifolds with bounded geometry. We
also prove existence of extremal under an extra condition. Finding
extremal of useful functionals is an useful problem in
mathematical analysis. For instance there is a vast literature
devoted to the study of ground state eigenvalues and
eigenfunctions which are extremal of the Dirichlet functional. The
Log Sobolev functional (\ref{logfanhan}) seems to be a mild
nonlinear perturbation of the Dirichlet functional. Indeed, they
share a common property i.e.
 there exist extremal functions for both functionals on compact domains or compact manifolds.
  However in the noncompact case
the similarity stops. For instance in $\R^n$, it is well known
that the Dirichlet functional does not have an extremal or $L^2$
eigenfunction. In contrast the Gaussian functions are extremals of
the Log Sobolev functional.  Over the years, Log Sobolev
inequality has found many applications in various branches of
mathematics and physics.  See for example the papers Gross
\cite{G:1}, \cite{G:2}, Federbush \cite{F:1},  Bakry and \'Emery
\cite{BE:1}, Bakry and Ledoux \cite{BL:1},  Diaconis and
Saloff-Coste \cite{DS:1} and Otto and Villani \cite{OV:1}.  A more
recent application was discovered by Perelman \cite{P:1} where he
introduced the fundamental $W$ entropy (\ref{Wshang}) and used it
as a key analytic tool to prove the Poincar\'e conjecture. The $W$
entropy is just the Log Sobolev functional (\ref{logfanhan})
scaled with certain time dependent parameter. For the Log Sobolev
functional, the
 existence problem of extremal functions in the compact case was solved by O. Rothaus \cite{Rot:1} 30 years ago.
 However, in the noncompact case, the problem is wide open.
 There has been no counter example or
 general existence result  for connected, noncompact manifolds with bounded geometry.
 We should mention that if one drops the connectedness,
  then it is easy to construct a manifold
 with infinitely many disconnected components, such that the Log Sobolev functional does not have an extremal.  See the example at the beginning of Section 3 e.g.  Also if the manifold is homogeneous
 such as $\R^n$, one can use symmetrization or translation or group action to
  prove existence of an extremal.

In addition to being an interesting problem in its own right, the
study of Log Sobolev inequality or $W$ entropy in the noncompact
setting is also important to Ricci flow. One reason is that many
of the more interesting singularity models are noncompact, even
when the Ricci flow under consideration is compact. One such
example in the three dimensional case is the round neck $S^2
\times R$, which is a typical singularity model. Using the
existence of extremals of his $W$ entropy, Perelman \cite{P:1}
proved a no breather theorem stating that shrinking breathers of
Ricci flows on compact manifolds are shrinking gradient solitons.
Recently, in the case $({\M}, g)$ is a noncompact
  gradient shrinking soltion, Carrillo and Ni \cite{CN:1} proved that
  potential functions are extremals for $W$ the entropy.

On p9, 3.2 Remark of the same paper, Perelman also wrote

{\it "Of course, this argument requires the existence of  minimizer,
and justifications of the integration by parts; this is easy if $M$ is closed, but can also be done with
more efforts on some complete $M$, ..."}

However, it is not known so far if the $W$ entropy always has an extremal for all noncompact manifolds
which are reasonably nice, such as those connected ones with bounded geometry.
The main theorem of the paper (Theorem \ref{thmain} or Theorem \ref{thWshang})  shows that on noncompact manifolds, the Log Sobolev functional or the $W$ entropy has an extremal function under a condition near infinity;
 it also shows
that an extremal function may not exist if the condition is violated, giving a negative answer to
the above question stipulated by Perelman. As another application
we partially extend Perelman's no breather theorem to the noncompact
case. See Section 4 below.

In order to state the result precisely, we first introduce a number of basic assumptions and notations.

{\bf Basic assumptions.} In this paper, unless stated otherwise,
we assume the $n$ dimensional Riemannian manifold ${\M}$ with
metric $g$ is a complete noncompact manifold with bounded geometry
which means:

1.  there exists a positive constant $\alpha$ such that
\[
| Rm | \le \alpha
\] where $Rm$ is the curvature tensor and $| Rm |$ is the maximum norm of $Rm$ under  $g$.

2. there exists a positive constant $\beta$ such that, for all $x
\in \M$,
\[
|B(x, 1)|_g \ge \beta.
\] Here $B(x, 1)$ is the geodesic ball of radius $1$, centered at $x$; and $|B(x, 1)|_g$ is the volume of
$B(x, 1)$ under the metric $g$.

It is well known that assumptions 1 and 2 imply that the
injectivity radius of $\M$ is bounded from below by a positive
constant.  See  \cite{CGT:1}  and \cite{CLY:1} e.g.

We will use the following notations throughout the paper.  $g_{ij}, R_{ij}$ will be the metric and Ricci
curvature;  $R$ is the scalar curvature; $\nabla$, $\Delta$ the corresponding gradient and
Laplace-Beltrami operator; $dg$ is the volume element; $c$, $C$ with or without index denote generic
positive constant that may change from line to line.

The main Log Sobolev inequality that we deal with in this paper is just the usual one perturbed by the scalar curvature of the manifold.  i.e.  there exist positive constant $a$ and another constant $c=c(a, {\M}, g)$ such that, for $v \in C^\infty_0(\M)$ and $\Vert v \Vert_{L^2(\M)}=1$,
\[
\int_{\M} v^2 \ln v^2 dg \le  a   \int_{\M} ( 4 | \nabla v|^2 + R v^2) dg + c(a, {\M}, g).
\]

The functional associated with the Log Sobolev inequality when $a=1$ is
\be
\lab{logfanhan}
L(v, g) \equiv \int_{\M} ( 4 | \nabla v|^2 + R v^2 - v^2 \ln v^2 ) dg,  \qquad \qquad v \in W^{1, 2}(\M).
\ee

One reason for involving the scalar curvature is,  after scaling
the functional by certain time dependent factor and coupled with
Ricci flow, it becomes Perelman's  $W$ entropy \cite{P:1},  which
is a fundamental quantity for Ricci flow. This relation is  shown
in (\ref{W=L}). The existence and nonexistence of extremal of the
Log Sobolev functional depends on two quantities given in the
definition below. The first one is just  the best Log Sobolev
constant or the infimum of the functional in (\ref{logfanhan}).
The second one is the best Log Sobolev constant at infinity. The
concept is motivated by P.L. Lions' concentration compactness
principle \cite{L:1}.

\begin{definition}
\lab{delamb} Let $({\M}, g)$ be a complete
noncompact manifold with bounded geometry.

 The best Log Sobolev constant of $({\M}, g)$ is the quantity
\[
\lambda=\lambda(\M)=\lambda({\M}, g) = \inf \{  \int_{\M} ( 4 |
\nabla v|^2 + R v^2 - v^2 \ln v^2 ) dg \quad | \, v \in
C^\infty_0(\M), \quad \Vert v \Vert_{L^2(\M)=1} \}.
\]

The best Log Sobolev constant of $({\M}, g)$ at infinity is the quantity
\[
\al
\lambda_\infty=\lambda_\infty({\M}, g)
= &\lim_{r \to \infty}
 \inf \{  \int_{\M-B(0, r)} ( 4 | \nabla v|^2 + R v^2 - v^2 \ln v^2 ) dg \quad | \\
 & \qquad \quad  \quad
v \in C^\infty_0(\M-B(0, r)), \quad \Vert v \Vert_{L^2(\M-B(0, r))=1} \}.
\eal
\]

Let $D$ be a domain in $\M$. The best Log Sobolev constant of $D$ is the quantity
\[
\lambda=\lambda(D)=\lambda(D, g) = \inf \{  \int_{D} ( 4 | \nabla
v|^2 + R v^2 - v^2 \ln v^2 ) dg \quad | \, v \in C^\infty_0(D),
\quad \Vert v \Vert_{L^2(D)=1} \}.
\]
\end{definition}

\begin{definition} (extremal)  Suppose $\lambda=\lambda({\M}, g)$ is a finite number. A function $v \in W^{1, 2}(\M)$ is called
an extremal  of the Log Sobolev functional  (\ref{logfanhan}) if $\Vert v \Vert_{L^2(\M)}=1$ and
\[
\int_{\M}  (4 | \nabla v|^2 + R v^2 - v^2 \ln v^2 ) dg = \lambda
\]
\end{definition}

The main result of the paper is the following theorem, or
equivalently Theorem \ref{thWshang} in Section 4.

\begin{theorem}
\lab{thmain}
  (a). Let $\M$ be a complete, connected
noncompact manifold with bounded geometry.
 Suppose $\lambda<\lambda_\infty$,
then there exists a smooth extremal $v$ for the Log Sobolev
functional in (\ref{logfanhan}). Also, there exist positive
constants $a, A>0$ and a point $0 \in \M$ such that
\[
v(x) \le A e^{- a d^2(x,0)}.
\]

(b).  There exists  a complete, connected
noncompact manifold with bounded geometry such that $\lambda = \lambda_\infty$, but
the Log Sobolev functional in (\ref{logfanhan}) does not have an extremal.
\end{theorem}

{\it Remark.}  Manifolds satisfying the condition $\lambda<\lambda_\infty$ are quite common.
For example, suppose $\M$ is asymptotically Euclidean, then
$
\lambda_\infty = \lambda({\R}^n)= \frac{n}{2} \ln(4 \pi) -n.
$ If there exists a compact domain $D \subset \M$ such that $\lambda(D)<\frac{n}{2} \ln(4 \pi) -n$.
Then
\[
\lambda \le \lambda(D)<\lambda_\infty.
\] It is easy to construct a domain such that $\lambda(D)$ is arbitrarily negative. One example
is the scaled  flat torus $h^2 (S^1 \times S^1) \times S^1$ when
the scaling factor $h \to 0$. See Lemma \ref{leLftube}.

Even though the Log Sobolev functional in the theorem contains the
scalar curvature $R$, the result still holds if one deletes the
scalar curvature. The proof only requires minor adjustment.

The rest of the paper is organized as follows.
Theorem \ref{thmain} (a) and (b) will be proven in Sections 2 and 3 respectively.
Applications on the $W$ entropy will be given in Section 4.

\section{Proof of Theorem \ref{thmain} (a), the existence part}

The proof of the theorem relies on the study of the
Euler-Langrange equation of the Log Sobolev functional:
 \be
\lab{eqln} 4 \Delta v - R v + 2 v \ln v + \lambda v =0. \ee When
$\lambda$ is the best Log Sobolev constant, this is the equation
satisfied by the extremal. Sometimes we also need to deal with
subsolutions to this equation.  A function $v \in W^{1,
2}_{loc}(\M)$ is called a subsolution to (\ref{eqln}) if it
satisfies the following inequality in the weak sense:
\be
\lab{subln}
 4 \Delta v - R v + 2 v \ln v + \lambda v \ge 0, \qquad \text{in} \qquad \M.
 \ee i.e., for any nonnegative, compactly supported test function $\psi$, we have
 \[
 \lambda \int_\M v \psi dg \ge  \int_\M  (4 \nabla v \nabla \psi + R v \psi - 2 \psi v \ln v  ) dg.
 \]

We will need a number of lemmas before proving the theorem.
The first lemma is a mean value type inequality for subsolutions of the above equation

\begin{lemma}
\lab{lemvi}
(a). Suppose $v$ is a bounded subsolution to the equation
(\ref{eqln}) in the ball $B(m, 2) \subset \M$ such that $\Vert v \Vert_{L^2(B(m, 2))} \le 1$. Here $ m \in \M$ which has bounded geometry. Then there exists a positive constant $C=C(n,
\alpha, \beta, \lambda)$ such that
\[
\sup_{B(m, 1)} v^2 \le C \int_{B(m, 2)} v^2 dg.
\]

(b).  Moreover, if $v$ is a bounded solution to (\ref{eqln}) in the ball $B(m, 2) \subset \M$ such that $\Vert v \Vert_{L^2(B(m, 2))} \le 1$, then
there exists a positive constant $C=C(n,
\alpha, \beta, \lambda)$ such that the gradient bound holds:
\[
\sup_{B(m, 1/2)} |\nabla v|^2 \le C \int_{B(m, 1)} v^2 dg.
\]

\proof
\end{lemma}

{\it Part (a).}
This part of the lemma and its proof is similar to that of Lemma 8.2.1 in \cite{Z:1} where the underlying manifold is
an $\epsilon$ horn. The proof relies on Moser's iteration and standard Sobolev inequality and takes advantage
of the slow growth of $\ln v$ when $v$ is large.

Given any $p \ge 1$, it is easy to see that
\begin{equation}
- 4 \Delta v^p + p R v^p \le 2 p v^p \ln v + p |\lambda| v^p.
\label{9.3}
\end{equation}

We select a smooth cut off function $\phi$ supported in $B(m, 2)$. Writing $w = v^p$ and using $w
\phi^2$ as a test function in (\ref{9.3}), we deduce
\[
4 \int \nabla (w \phi^2) \nabla w dg + p \int R (w \phi)^2 dg \le 2 p \int (w \phi)^2 \ln v dg +p \int |\lambda|
(w \phi)^2 dg.
\] By the bound on the curvature tensor $|Rm| \le \alpha$, we deduce
\[
4 \int \nabla (w \phi^2) \nabla w dg \le  p \int (w \phi)^2 \ln v^2 dg + (C \alpha +|\lambda|) p \int (w \phi)^2 dg,
\]which induces, after integration by parts,
\begin{equation}
4 \int |\nabla (w \phi)|^2 dg \le  4 \int |\nabla
\phi|^2 w^2 dg  + (C \alpha p +|\lambda|) \int (w \phi)^2 dg + p \int (w \phi)^2 \ln v^2 dg.
\label{9.4}
\end{equation}

We need to dominate the  last term in (\ref{9.4}) by the left hand
side of (\ref{9.4}). For one positive number $a$ to be chosen later,
it is clear that
\[
\ln v^2 \le v^{2a} + c(a).
\]Hence for any fixed $q>n/2$, the H\"older inequality implies
\[
\begin{aligned}
p \int (w \phi)^2 \ln v^2 dg &\le p \int (w \phi)^2 v^{2a} dg+ p c(a) \int (w \phi)^2 dg\\
&\le p \left( \int v^{2 a q} dg \right)^{1/q} \   \left( \int (w \phi)^{2  q/(q-1)} dg
\right)^{(q-1)/q}
+ p c(a) \int (w \phi)^2 dg.
\end{aligned}
\]We take $a = 1/q$ so that $2 a q =2$.
Since the $L^2$ norm of $u$ is less than $1$ by assumption, the above implies
\[
p \int (w \phi)^2 \ln v^2 dg \le  p  \left( \int (w \phi)^{2  q/(q-1)} dg
\right)^{(q-1)/q}
+ p c(a) \int (w \phi)^2 dg.
\]By interpolation inequality (see p84  \cite{HL:1} e.g.), it holds, for any $b>0$,
\[
\left( \int (w \phi)^{2  q/(q-1)} dg \right)^{(q-1)/q}
\le b
\left( \int (w \phi)^{2 n/(n-2)} dg \right)^{(n-2)/n}  + c(n, q) b^{-n/(2q-n)}
 \int (w \phi)^2 dg.
\]Therefore
\begin{equation}
 p \int (w \phi)^2 \ln v^2 dg  \le  p b \left( \int (w \phi)^{2
n/(n-2)} dg \right)^{(n-2)/n}  + C p b^{-n/(2q-n)} \int (w
\phi)^2 dg + p C \int (w \phi)^2 dg.
\label{wphiln}
\end{equation}

Since the manifold $\M$ has bounded geometry, it is well known
(\cite{Au:1},\cite{Heb:1}, \cite{HV:1} and \cite{Sa:1} e.g.) that
a standard Sobolev inequality holds, i.e. there exist positive
constants $S_0$ depending on $\alpha, \beta, n$ such that
\[
S_0 \left(\int (w  \phi)^{2n/(n-2)} dg \right)^{(n-2)/n} \le  \int |\nabla (w \phi)|^2 dg +
 \int  (w \phi)^2 dg.
 \] This and (\ref{9.4}) imply
\begin{equation}
S_0 \left(\int (w  \phi)^{2n/(n-2)} dg \right)^{(n-2)/n} \le   \int |\nabla
\phi|^2 w^2 dg  + (C \alpha p +|\lambda|+1) \int (w \phi)^2 dg + p \int (w \phi)^2 \ln v^2 dg.
\label{9.5}
\end{equation}

Substituting (\ref{wphiln}) to the right hand side of
(\ref{9.5}), we arrive at
\[
\al
 S_0 \left(\int (w \phi)^{2n/(n-2)} dg \right)^{(n-2)/n} & \le 4 \int |\nabla \phi|^2 w^2 dg
  +
 p b \left( \int (w \phi)^{2 n/(n-2)} dg \right)^{(n-2)/n} \\
 &\qquad + c(n, q) p b^{-n/(2q-n)}
 \int (w \phi)^2 dg
+ p c(a) \int (w \phi)^2 dg.
\eal
\]Take $b$ so that $p b = S_0/2$. It is clear that exist positive constant
$c =c(S_0, n, q)$
and $p_0=p_0(n, q)$ such that
\begin{equation}
\left(\int (w \phi)^{2n/(n-2)} dg \right)^{(n-2)/n} \le c
(p+1)^{p_0} \int (|\nabla \phi|^2+1) w^2 dg.
\label{9.7}
\end{equation}

From here one can use standard Moser's iteration to prove the desired bound.
We briefly sketch the main steps.
Let $\xi_k=\xi_k(s)$, $k=0, 1, 2, ...$, be a smooth one variable function such that
$\xi_k(s)=1$ when $s \in [0, 1+(1/2^{k+1})]$; $0 \le \xi_k(s) \le 1$,  when $s \in [1+(1/2^{k+1}),
1+ (1/2^k) ]$
and $\xi_k(s) =0$,  when $s \in [1+(1/2^k), \infty)$. We also require that
$|\xi'(s)| \le c/2^k$.
Set the test function $\phi_k=\xi_k(d(x, m))$. Then it is clear that
\begin{equation}
|\nabla \phi_k| \le \frac{c}{2^k}. \label{9.8}
\end{equation}
By (\ref{9.7}) and (\ref{9.8})
\begin{equation}
\left(\int_{B(m, 1+ (1/2^{k+1}))} w^{2n/(n-2)} dg \right)^{(n-2)/n} \le \frac{C}{2^{2k}}
 (p+1)^{p_0} \int_{B(m, 1+(1/2^k))} w^2 dg.
\label{9.9}
\end{equation}

Recall that $w = v^p$. We iterate (\ref{9.9}) with $p=(n/(n-2))^k$,
$k=0, 1, 2, ...$ Following Moser, we get
\[
\sup_{B(m, 1)} v^2 \le C  \int_{B(m, 2)} v^2 dg.
\] This proves part (a) of the lemma. \\

{\it Part (b).}  By standard computation, in local orthonormal system, we have
\[
\Delta |\nabla v|^2 = 2 \Sigma_{i, j} v^2_{ij} + 2 \Sigma_i ( \Delta v)_i v_i +
4 R_{ij} v_i v_j.
\] Here $v_i$ is the covariant derivative of $v$ and $R_{ij}$ is the Ricci curvature.
Since $v$ is a solution to (\ref{eqln}), we know that
\[
 ( \Delta v)_i  v_i= \frac{1}{4} ( R v - 2 v \ln v - \lambda v)_i v_i =
  \frac{1}{4} ( R_i v v_i + R v^2_i - 2 v^2_i \ln v -  2 v^2_i - \lambda v^2_i).
\] Since , by part (a), $v \le C$ in $B(m, 1)$ , we have $-\ln v \ge -\ln C$. Hence there exists a positive constant $C$ such that
\[
\Delta |\nabla v|^2 \ge -C (|\nabla v|^2 + v^2)
\] in the ball $B(m, 1)$. From here,  we can use Moser's iteration for standard Laplacian to conclude
that
\[
\sup_{B(m, 1/2)} |\nabla v|^2 \le C \int_{B(m, 2r/3)} (|\nabla v|^2 + v^2) dg \le C  \int_{B(m, r)}  v^2 dg.
\]
\qed

The next lemma shows that interior maximum value of a positive solution of equation (\ref{eqln}) in a ball
has a positive lower bound independent of the ball.  This property in case of compact manifolds
was already observed in Section 17.2 of \cite{CCGGIIKLLN3:1}.

\begin{lemma}
\lab{lexiajie}
Let $v$ be a smooth positive solution of equation (\ref{eqln}) in the ball $B(0, r) \subset \M$ such that
$v=0$ on $\partial B(0, r)$. Here
$0$ is a point in $\M$ and $r>0$. Then
\[
\sup_{B(0, r)} v \ge e^{(\inf R-\lambda)/2}.
\] i.e. the maximum value of $v$ is bounded from below by a positive constant depending only
on $\lambda$ and the lower bound of the scalar curvature.
\proof
\end{lemma}

Since $v$ is $0$ at the boundary, clearly the maximum of $v$ is reached at some point $x_0$ in the interior of the ball $B(0, r)$.
Hence $\Delta v(x_0) \le 0$, which implies, by equation (\ref{eqln}),
\[
- R(x_0) v(x_0) + 2 v(x_0) \ln v(x_0) + \lambda v(x_0) \ge 0.
\] From this, the lemma follows. \qed

\begin{lemma}
\lab{ledecay}
 Let $v$ be a bounded  subsolution to (\ref{eqln}) on $\M$ such that
$\Vert v \Vert_{L^2(\M)} \le 1$. Let $0$ be a reference point on $\M$.  Then there exist
positive numbers $r_0$, $a$ and $A$, which may depend on $\alpha, \beta$ and the location of the reference point such that
\[
v(x) \le A e^{-a d^2(x, 0)}, \qquad \text{when} \qquad d(x, 0) \ge r_0.
\]
\proof
\end{lemma}

Recall from Lemma \ref{lemvi} that there exists a constant $C>0$ such that
\[
v^2(x) \le C \int_{B(x, 2)} v^2 dg, \qquad x \in \M.
\]This infers
\[
- 2 \ln v(x) \ge - \ln C - \ln \int_{B(x, 2)} v^2 dg.
\]Since $\int_\M v^2 dg \le 1$, we know that
\[
\lim_{d(x, 0) \to \infty} \int_{B(x, 1)} v^2 dg =0.
\] Therefore $-\ln v(x) \to +\infty$ when $d(x, 0) \to \infty$. Thus, there exists $r_0>0$, such that,
when $d(x, 0) \ge r_0$, we have
\be
R(x) - \ln v(x) - \lambda \ge 0, \qquad \text{and} \quad  v(x) \le e^{-1}
\ee Substituting this to (\ref{subln}), we deduce,
\[
4 \Delta v(x) + v(x) \ln v(x) \ge v(x) (R(x)  -\ln  v(x) -\lambda)
\ge 0.
\] Hence,
when $d(x, 0) \ge r_0$, we have
\be
\lab{ddv>}
4 \Delta v(x) + v(x) \ln v(x) \ge 0, \qquad  \text{and} \quad v(x) \le e^{-1}.
\ee

Next we compare $v$ with a model function
\be
\lab{Jx}
J=J(x)= e^{- a L^2(x) + a r^2_0 -1}.
\ee Here $a>0$ is to be decided later; $L=L(x)$ is a smooth function on $\M$, which satisfies
\[
|\nabla L(x) | \le C_1, \qquad |\nabla^2 L(x) | \le C_1, \qquad x \in \M,
\]
\[
C^{-1}_1 L(x) \le d(x, 0) \le C_1 L(x), \qquad d(x, 0) \ge r_0.
\] Under our assumption of bounded geometry, it is well known that such a function exists.
For instance, let $\eta \ge 0$ be a smooth function in
$C^\infty_0({\R}^n)$, supported in a ball centered at the origin,
whose radius is less than the injectivity radius of $\M$. If also
$\Vert \eta \Vert_{L^1({\R}^n)}=1$, then \be \lab{Lx=} L(x) =
\int_{{\R}^n} \eta(w) [ d(0, exp_x(w)) + 1] dw \ee satisfies the
above requirements. See also the proof of Proposition 19.37 in
\cite{CCGGIIKLLN3:1}, e.g. Since $d(x, 0)$ and $L(x)$ are
comparable when they are large, by (\ref{ddv>}), we can choose
$r_0$ sufficiently large so that \be \lab{ddv>L} 4 \Delta v(x) +
v(x) \ln v(x) \ge 0, \qquad  \text{and} \quad v(x) \le e^{-1} \ee
when $L(x) \ge r_0$.

By direct computation
\[
\Delta J = J [ 4 a^2 | \nabla L |^2 L^2 - 2 a L  \Delta L - 2 a | \nabla L|^2 ],
\]
\[
J \ln J = J ( - a L^2 + a r^2_0 -1).
\]Hence
\[
\al
4 \Delta J + J \ln J &= J  [ 16 a^2 | \nabla L |^2 L^2 - 8 a L \Delta L - 8 a | \nabla L|^2
- a L^2 + a r^2_0 -1 ]\\
&\le J  [ 16 a^2 C^2_1 L^2 + 8 a C_1 L
- a L^2 + a r^2_0 -1 ].
\eal
\]This implies, for some $C_2>0$,
\[
4 \Delta J + J \ln J \le  J  [ C_2 a^2 L^2
- a L^2 + a r^2_0 - (1/2) ].
\]
We take $a = \min \{ \frac{1}{C_2}, \frac{1}{\sqrt{2 C_2 r^2_0}} \}$. Then
\[
4 \Delta J + J \ln J \le 0
\]when $L(x) \ge r_0$ and $J(x)= e^{-1}$ when $L(x)=r_0$.
This and (\ref{ddv>L}) show that
\[
\al
\begin{cases}
4 \Delta (J-v) + J \ln J - v \ln v  \le 0, &\text{if} \qquad L(x) \ge r_0,\\
J(x) \le e^{-1}, \qquad v(x) \le e^{-1},  &\text{if} \qquad L(x) \ge r_0\\
(J-v)(x) \ge 0, &\text{if} \qquad L(x)=r_0,\\
(J-v)(x) \to 0,  &\text{if} \qquad L(x) \to \infty,\\
\end{cases}
\eal
\]Since $J(x), v(x) \le e^{-1}$, by the mean value theorem, there exists a function $f=f(J(x), v(x))$,
$0 < f \le e^{-1}$ such that
\[
J(x) \ln J(x) - v(x) \ln v(x) = (\ln f + 1) (J(x)-v(x)).
\]Observe that
\[
\ln f + 1 \le ln e^{-1} + 1 \le 0, \qquad \text{when} \quad L(x)
\ge r_0.
\]Therefore we can apply the standard maximum principle for the elliptic
inequality on
\[
4 \Delta (J-v)(x) + (\ln f + 1) (J-v)(x) \le 0, \qquad \text{when}
\quad  L(x) \ge r_0
\]to conclude that
\[
v(x) \le J(x) = e^{-a L^2(x) + a r^2_0 -1}, \qquad \text{when}
\quad L(x) \ge r_0.
\]Since $L(x)$ and $d(x, 0)$ are comparable when they are large, we have proven the lemma by
making $a$ smaller if necessary.
\qed

\begin{lemma}
\label{leLMD} Let $({\bf M}, g)$ be a complete noncompact manifold with bounded geometry.
Let $v \in W^{1, 2}(\M)$, $\Vert v \Vert_{L^2(\M)}=1$ be  a bounded sub-solution of (\ref{eqln}) i.e.
 \[
 4 \Delta v - R v + 2 v \ln v + \lambda v \ge 0.
 \] Here $\lambda$ is a  constant. Let $D$ be
a bounded domain in $\M$
 and define
\begin{equation}
\lambda (D) = \inf \{ \int ( 4 |\nabla v|^2 + R v^2- v^2 \ln v^2) dg \ |
\ v \in C^\infty_0(D), \ \Vert v \Vert_2 = 1 \},
\label{9.12}
\end{equation}
 For any smooth cut-off
function $\eta \in C^\infty_0(D)$, \ $0 \le \eta \le 1$, it
holds
\[
\lambda (D) \int (v \eta)^2 dg \le \lambda \int (v \eta)^2 dg+ 4 \int  v^2 |\nabla \eta|^2 dg -
\int (v \eta)^2 \ln \eta^2 dg.
\]
\proof
\end{lemma}

Since $\eta v/\Vert \eta v \Vert_2 \in C^\infty_0(D)$  and its $L^2$ norm is $1$,
we have, by  definition,
\[
\lambda(D) \le \int \left[ 4 \frac{|\nabla (\eta v)|^2 }{\Vert \eta v \Vert^2_2 } +
R \frac{(\eta v)^2}{\Vert \eta v \Vert^2_2 } - \frac{(\eta v)^2}{\Vert \eta v \Vert^2_2 }
\ln \frac{(\eta v)^2}{\Vert \eta v \Vert^2_2 } \right] dg.
\]This implies
\begin{equation}
\lambda(D) \Vert \eta v \Vert^2_2  \le \int \left[ 4 |\nabla (\eta
v)|^2 +R (\eta v)^2 -
 (\eta v)^2 \ln (\eta v)^2 \right] dg  + \Vert \eta v \Vert^2_2 \ln  \Vert \eta v \Vert^2_2.
 \label{9.14}
\end{equation}

 On the other hand, $v$ satisfies
 \[
 4 \Delta v - R v + 2 v \ln v + \lambda v \ge 0.
 \]Using $\eta^2 v$ as a test function here, we deduce
 \[
 \lambda  \int (\eta v)^2 dg \ge - 4 \int (\Delta v) \eta^2 v dg + \int R (\eta v)^2 dg -
 2 \int (\eta v)^2 \ln v dg.
 \]By direct calculation
 \[
 -4 \int (\Delta v)  \eta^2 v dg= 4 \int |\nabla (\eta v)|^ 2 dg- 4 \int v^2 |\nabla \eta|^2 dg.
 \]Hence
 \begin{equation}
  \lambda \int (\eta v)^2 dg \ge 4 \int |\nabla (\eta v)|^ 2 dg- 4 \int v^2 |\nabla \eta|^2 dg
 + \int R (\eta v)^2 dg-
 2 \int (\eta v)^2 \ln v dg.
 \label{9.15}
 \end{equation} Comparing (\ref{9.15}) with (\ref{9.14}) and noting that $\Vert \eta v \Vert_2 < 1$, we
 obtain
 \[
 \lambda(D) \Vert \eta v \Vert^2_2  \le  \lambda \Vert \eta v \Vert^2_2 +
 4 \int |\nabla \eta|^2 v^2 dg - \int (\eta v)^2 \ln \eta^2 dg.
 \]
 \qed

 The next lemma is a stability result for the infimum of the Log Sobolev functional under $C^2$
 perturbation of the metric. We believe it should be known. However, since we can not find it in the
 literature, we present it here.

 \begin{lemma}
 \lab{leLstab}
 Let $D \subset \M$ be a compact domain. For any $\epsilon>0$, there exists $\delta>0$ such that the following
 statement is true.

 Let $g_1$ and $g_2$ be two metrics on $\M$ such that
 \[
 \Vert g_1 - g_2 \Vert_{C^2(D, g_1)} <\delta.
 \]Here $\Vert \cdot \Vert_{C^2(D, g_1)}$ stands for the $C^2$ norm for $(2, 0)$ tensor fields under the
 metric $g_1$, restricted to the domain $D$. Then
 \[
 |\lambda(D, g_1) - \lambda(D, g_2)|<\epsilon.
 \] Here, for $i=1, 2$,
 \[
 \lambda(D, g_i) =  \inf \{  \int_{D} ( 4 | \nabla_{g_i} v|^2 + R_{g_i} v^2 - v^2 \ln v^2 ) dg_i \quad | \quad
v \in C^\infty_0(D), \quad \Vert v \Vert_{L^2(D, g_i)=1} \}.
\]
\proof
\end{lemma}

By definition of $ \lambda(D, g_1) $, there exists a function $v \in C^\infty_0(D)$ such that
$\Vert v \Vert_{L^2(D, g_1)} = 1$ and that
\[
\lambda(D, g_1) + \epsilon >  \int_{D} ( 4 | \nabla_{g_1} v|^2 + R_{g_1} v^2 - v^2 \ln v^2 ) dg_1.
\] Recall, in  local coordinate patch $U$ with coordinate $\{ x^1, ..., x^n \}$,
\[
| \nabla_{g_1} v|^2 = g^{ij}_1 \partial_i v \partial_j v.
\] Hence, in each local coordinate patch,
\[
-\epsilon< | \nabla_{g_1} v|^2 - | \nabla_{g_2} v|^2 < \epsilon; \qquad |R_{g_1}- R_{g_2}|<\epsilon;
\qquad |dg_1-dg_2|<\epsilon
\]when $\Vert g_1 - g_2 \Vert_{C^2(D, g_1)} <\delta$ with $\delta$ being sufficiently small. Since $D$ is compact, it can be covered by finitely many local
charts. Therefore, there exists $C>0$ such that
\[
\lambda(D, g_1) + \epsilon >  \int_{D} ( 4 | \nabla_{g_2} v|^2 + R_{g_2} v^2 - v^2 \ln v^2 ) dg_2 - C \epsilon.
\] Consider the function $\tilde v = v/\Vert v \Vert_{L^2(D, g_2)}$.
Then the above inequality becomes
\[
\lambda(D, g_1) + \epsilon >  \int_{D} ( 4 | \nabla_{g_2} \tilde v|^2 + R_{g_2} \tilde v^2 - \tilde v^2 \ln \tilde v^2 ) dg_2  \, \Vert v \Vert^2_{L^2(D, g_2)} - \Vert v \Vert^2_{L^2(D, g_2)} \ln \Vert v \Vert^2_{L^2(D, g_2)}   - C \epsilon.
\] Since $\Vert \tilde v \Vert_{L^2(D, g_2)} =1$, we deduce
\[
\lambda(D, g_1) + \epsilon > \lambda(D, g_2)  \, \Vert v \Vert^2_{L^2(D, g_2)} - \Vert v \Vert^2_{L^2(D, g_2)} \ln \Vert v \Vert^2_{L^2(D, g_2)}   - C \epsilon.
\] Notice that $  \Vert v \Vert^2_{L^2(D, g_1)} =1$ and $\Vert g_1 - g_2 \Vert_{C^2(D, g_1)} <\delta$.
Thus $| 1 -  \Vert v \Vert^2_{L^2(D, g_2)} |< \epsilon$ when $\delta$ is sufficiently small.
Hence there exists $C>0$ such that
\[
\lambda(D, g_1) + C \epsilon > \lambda(D, g_2).
\] In the same manner, we obtain
\[
\lambda(D, g_2) + C \epsilon > \lambda(D, g_1)
\]  which shows
\[
| \lambda(D, g_1)  - \lambda(D, g_2) | < C \epsilon.
\]
\qed

Now we are ready to give the
\medskip

{\bf Proof of Theorem \ref{thmain} (a), the existence part.}

We assume $\lambda < \lambda_\infty$.  First we prove that $\lambda$ is finite. Since $\M$ has bounded
geometry, it is well known (c.f. \cite{Au:1},  \cite {Heb:1}, \cite{HV:1}) that the following Sobolev inequality holds:
there exist positive constants $S_0$
depending on $\alpha, \beta, n$ such that, for all $v \in C^\infty_0(\M)$,
\[
S_0 \left(\int v^{2n/(n-2)} dg \right)^{(n-2)/n} \le  \int |\nabla v |^2 dg +
 \int  v^2 dg.
\] Under the assumption $\Vert v \Vert_{L^2(\M)} =1$, a quick application of Jensen's inequality on
the Sobolev inequality shows, for a constant $C=C(n, S_0)$ and all $\epsilon>0$,
\[
\int v^2 \ln v^2 dg \le \epsilon^2 \int | \nabla v |^2 dg - \frac{n}{2} \ln \epsilon^2 + \epsilon^2 + C.
\]Taking  $\epsilon =2$ and using the assumption that the scalar curvature $R$ is bounded, we deduce
\be
\lab{l>-wq}
\lambda =  \inf \{ \int_{\M} ( 4 | \nabla v|^2 + R v^2 - v^2 \ln v^2 ) dg \, | \,
v \in C^\infty_0(\M), \, \Vert v \Vert_{L^2(\M)} =1 \} > -\infty,
\ee i.e. $\lambda$  is finite.

For positive integers $k$, consider the domains
\[
D(0, k)= \{ x \in \M \, | \, L(x)<k \}
\]where $L=L(x)$ is the smooth function defined by (\ref{Lx=}), which is comparable to $d(0, x)$ when
it is large.  By properties of $L=L(x)$, $\partial D$ is a $C^2$ boundary.
Given a positive integer $k$, let $\lambda_k$ be the best Log Sobolev constant of the ball $D(0, k)$, i.e.
\[
\lambda_k = \lambda(D(0, k)) =  \inf \{ \int ( 4 |\nabla v|^2 + R v^2- v^2 \ln v^2) dg \ |
\ v \in C^\infty_0(D(0, k)), \ \Vert v \Vert_2 = 1 \}.
\] According to \cite{Rot:1}, $\lambda_k$ is finite and there exists a smooth extremal function $v_k$ on
$D(0, k)$, which satisfies
\[
\al
\begin{cases}
4 \Delta v_k - R v_k + 2 v_k \ln v_k + \lambda_k v_k =0, \qquad \text{in} \qquad D(0, k)\\
v_k=0, \qquad \text{on} \qquad \partial D(0, k).
\end{cases}
\eal
 \]
 We  mention that $v_k$ is uniformly bounded in $C^\alpha(\M)$ norm, i.e., there exists a
 positive constant $C$ such that
 \be
 \lab{vkwinfty}
 \Vert v_k \Vert_{C^{\alpha}( D(0, k))} \le C.
 \ee A proof goes as follows.
 We extend $v_k$ to a function on the whole manifold $\M$ by setting $v_k(x)=0$ when $x
 \in \M - D(0, k)$. The extended function is still denoted by $v_k$. Then $v_k \in W^{1, 2}(\M)$, and $v_k$ satisfies the following inequality in the weak sense
 \[
 4 \Delta v_k - R v_k + 2 v_k \ln v_k + \lambda_k v_k \ge 0, \qquad \text{in} \qquad \M.
 \]i.e., for any nonnegative, compactly supported test function $\psi$, we have
 \[
 \lambda_k \int_\M v_k \psi dg \ge  \int_\M  (4 \nabla v_k \nabla \psi + R v_k \psi - 2 \psi v_k \ln v_k  ) dg.
 \] By Lemma \ref{lemvi},  the norm
 $\Vert v_k \Vert_{L^\infty(\M)}$ is uniformly bounded. Hence the original $v_k$ in $D(0, k)$ is actually a bounded weak solution
 to the Poisson equation
 \[
 \begin{cases}
 \Delta v_k(x) = f_k(x), \qquad x \in D(0, k) \\
 v_k(x)=0, \quad x \in \partial D(0, k)
 \end{cases}
 \] with $\Vert f_k \Vert_{L^\infty(\M)} \le C$. Note that $\partial D(0, k)$ is given by $L(x)=k$ and
 $|\nabla L(x)| + |\nabla^2 L(x)| \le C$ when $L(x)$ is large. Thus  $\partial D(0, k)$ is $C^2$ boundary which can be expressed by a uniform $C^2$ function locally in geodesic balls
 of radius less than the injectivity radius of $\M$. Hence the standard elliptic theory
 shows (\ref{vkwinfty}) is true.

 By (\ref{l>-wq}), $\lambda_k \ge \lambda>-\infty$ and $\{ \lambda_k \}$ is a decreasing sequence.
 Hence $\{ \lambda_k \}$ is uniformly bounded by a number, say $\Lambda$.
 According to Lemma \ref{lexiajie}, there exists a point $x_k \in D(0, k)$ and a uniform constant
 $C=C(n, \alpha, \beta, \Lambda)>0$ such that
 \be
 \lab{vxk>c}
 v_k(x_k) \ge C>0, \qquad k=1, 2, ...
 \ee

We consider $2$ cases. \\

{\it Case 1.}  $\{ x_k \}$ is a bounded sequence in $\M$, i.e. $d(x_k, 0)$ is uniformly bounded. \\

By Lemma \ref{lemvi},  the sequence $\{ v_k \}$ of extended functions is uniformly bounded in $L^\infty$ norm, $k=2, 3, ...$. By (\ref{vkwinfty}) we can find a subsequence, still denoted by $\{ v_k \}$, which converges in $C^\alpha_{loc}$ norm to a smooth, nonnegative function $v \in C^\infty(\M)$ that solves the equation
\[
4 \Delta v - R v + 2 v \ln v + \lambda v = 0.
\] The lower bound in (\ref{vxk>c}) ensures that $v$ is a positive solution. Moreover $\Vert v \Vert_{L^2(\M)} \le 1 $ by Fatou's Lemma. By Lemma \ref{ledecay}, there exist positive constants $a $ and $A$ such that
\[
v(x) \le A e^{- a d^2(x, 0)} \qquad x \in \M.
\]The classical volume comparison theorem tells us that $|B(0, k)|_g$ grows at most like $e^{c k}$, where
$c$ depends on the curvature bound $\alpha$ and $n$.  Hence we can multiply the above equation by $v$ and perform integration by parts to deduce
\be
\lab{lvg=lv2}
L(v, g) = \int_{\M}  [4 |\nabla v|^2 + R v^2  -v^2 \ln v^2 ] dg = \lambda \int_{\M} v^2 dg.
\ee

If $\int_{\M} v^2 dg = 1$, then $v$ is an extremal function of the
Log Sobolev functional $L$, and the proof Theorem \ref{thmain} (a)
is done. So we suppose $\int_{\M} v^2 dg < 1$.  We consider the
function
\[
\tilde v = \frac{v}{\Vert v \Vert_{L^2(\M)}}.
\]Then $\Vert \tilde v \Vert_{L^2(\M)} =1$ and (\ref{lvg=lv2}) infers
\[
\al
\lambda &= L(v, g) \Vert v \Vert^{-2}_{L^2(\M)} = \frac{\int_{\M}  [4 |\nabla v|^2 + R v^2  -v^2 \ln v^2 ] dg }{
\Vert v \Vert^{2}_{L^2(\M)}}\\
&=\int_{\M}  [4 |\nabla \tilde v|^2 + R \tilde{v}^2  -\tilde{v}^2 \ln \tilde{v}^2 ] dg - \ln  \Vert v \Vert^{2}_{L^2(\M)} \\
& \ge \lambda - \ln  \Vert v \Vert^{2}_{L^2(\M)} .
\eal
\] The last step is due to the definition that $\lambda$ is the infimum of the Log Sobolev functional.
If the assumption  $\int_{\M} v^2 dg < 1$ is valid, we would get the contradiction $\lambda>\lambda$.
Hence $\int_{\M} v^2 dg = 1$ and $v$ is indeed an extremal. This finishes the proof in Case 1.

\medskip
{\it Case 2.}  $\{ x_k \}$ is an unbounded sequence in $\M$.
\medskip

Since $\M$ has bounded geometry,  by Hamilton's  compactness theorem, the pointed manifolds $(\M, x_k, g)$ converges in $C^\infty_{loc}$
topology (also called Cheeger-Gromov sense), to a complete limit manifold $(M_\infty, x_\infty, g_\infty)$.
This limit manifold also has bounded geometry.

Recall $v_k (\ge 0)$ solves
\[
\al
\begin{cases}
4 \Delta v_k - R v_k + 2 v_k \ln v_k + \lambda_k v_k =0, \qquad \text{in} \qquad D(0, k)\\
v_k=0, \qquad \text{on} \qquad \partial D(0, k).
\end{cases}
\eal
 \]  We extend $v_k$ to a function on the whole manifold $\M$ by setting $v_k(x)=0$ when $x
 \in \M - D(0, k)$. The extended function is still denoted by $v_k$. Then, as in Case 1,
 $v_k \in C^\alpha(\M)
 \cap W^{1, 2}(\M)$, and $v_k$ satisfies the following inequality in the weak sense
 \[
 4 \Delta v_k - R v_k + 2 v_k \ln v_k + \lambda_k v_k \ge 0, \qquad \text{in} \qquad \M.
 \] Since $v_k$ is nonnegative and uniformly bounded by Lemma \ref{lemvi}, the standard elliptic theory
 shows that a subsequence of $\{ v_k \}$, converges in $C^\alpha_{loc}$ sense to a function
 $v \in C^\alpha(\M_\infty)
 \cap W^{1, 2}(\M_\infty)$.  Moreover $v$
satisfies the following inequality in the weak sense
 \[
 4 \Delta v - R v + 2 v \ln v + \lambda v \ge 0, \qquad \text{in} \qquad \M_\infty.
 \]i.e., for any nonnegative, compactly supported test function $\psi$, we have
 \[
 \lambda \int_{\M_\infty} v \psi dg_\infty \ge
  \int_{\M_\infty}  (4 \nabla v \nabla \psi + R v \psi - 2 \psi v \ln v  ) dg_\infty.
 \] Here the Laplacian $\Delta$, the gradient $\nabla$ and the scalar curvature $R$ are with respect to
 the limiting metric $g_\infty$.  Since $v_k(x_k)$ converges to $v(x_\infty)$, by (\ref{vxk>c}), we also know that
 \be
 \lab{vxinfty}
 v(x_\infty)> C >0.
 \ee By Lemma \ref{ledecay} and Fatou Lemma,  there hold the bounds
 \be
 \lab{vinftydecay}
 v(x) \le A e^{- a d^2(x, x_\infty, g_\infty)}, \quad x \in \M_\infty; \qquad \int_{\M_\infty} v^2(x) dg_\infty \le 1.
 \ee

 Let $r>0$ be a large number to be fixed later. Define, on the manifold $(M_\infty, g_\infty)$ and under the metric $g_\infty$,
\[
 \lambda (B(x_\infty, r)) = \inf \{ \int ( 4 |\nabla v|^2 + R v^2- v^2 \ln v^2) dg_\infty \ |
\ v \in C^\infty_0(B(x_\infty, r)), \ \Vert v \Vert_2 = 1 \}.
\]
 We choose a smooth cut-off
function $\eta \in C^\infty_0(B(x_\infty, r))$ such that \ $0 \le \eta \le 1$, $\eta =1 $ on $B(x_\infty, r/2)$ and that $|\nabla \eta| \le C/r$.
By Lemma \ref{leLMD}, it
holds
\be
\lab{lamBxinfty}
\lambda (B(x_\infty, r)) \le \lambda + 4 \frac{ \int  v^2 |\nabla \eta|^2 dg_\infty}{\int (v \eta)^2 dg_\infty} -
\frac{\int (v \eta)^2 \ln \eta^2 dg_\infty}{\int (v \eta)^2 dg_\infty}.
\ee By (\ref{vxinfty}) and the fact that $v$ is in $C^\alpha(M_\infty)$, we can find a positive constant
$c>0$ such that
\[
\int (v \eta)^2 dg_\infty \ge \int_{B(x_\infty, r/2)} v^2 dg_\infty \ge c.
\] From this and (\ref{lamBxinfty}), using properties of $\eta$, we deduce
\[
\lambda (B(x_\infty, r)) \le \lambda + C(1 + 1/r) \int_{B(x_\infty, r)-B(x_\infty, r/2)} v^2 dg_\infty.
\] By (\ref{vinftydecay}) and the classical volume comparison theorem, this implies
\[
\lambda (B(x_\infty, r)) \le \lambda + C(1 + 1/r) e^{-a r^2/4} e^{c \alpha r}.
\] Here, as before $\alpha$ is the bound on the curvature tensor.
Thus, for any $\epsilon>0$, there exists $r_0>0$ such that
\be
\lab{lamb>lambB}
\lambda =\lambda(\M) \ge \lambda (B(x_\infty, r))  -\epsilon
\ee when $r \ge r_0$.

By definition of $(M_\infty, x_\infty, g_\infty)$ as a limit manifold, for any $\delta>0$, when $k$ is sufficiently large, there
exists a diffeomorphism $F$  from $B(x_\infty,  r)$ onto an open set $U \subset \M$, which contains $x_k$, such that $(F^*)^{-1} g_\infty$ and $g$ are $\delta$ close in $C^\infty$ topology, when they are restricted to $U$. By Lemma \ref{leLstab}, we have, when $\delta$ is sufficiently small,
\be
\lab{lamb>lambU}
\lambda(B(x_\infty,  r)) = \lambda(B(x_\infty,  r), g_\infty) = \lambda(U, (F^*)^{-1} g_\infty) >\lambda(U,  g) -  \epsilon.
\ee  By definition of $U$, we know that  for any $x \in U$,
\[
d(x, x_k, (F^*)^{-1} g_\infty)< r
\] which implies, since $(F^*)^{-1} g_\infty$ and $g$ are $\delta$ close,
\[
d(x, x_k,  g)< (1+ C \sqrt{\delta}) r.
\]Hence, when $\delta$ is sufficiently small, it holds
\[
U \subset B(x_k, 2r, g).
\] This and (\ref{lamb>lambU}) tell us that
\[
\lambda(B(x_\infty,  r))  > \lambda(B(x_k, 2r, g),  g) -  \epsilon.
\]Recall that $d(x_k, 0, g) \to \infty$ when $k \to \infty$. Therefore, when $k$ is large,
\[
B(x_k, 2r, g) \subset \M  - B(0, d(x_k, 0, g)/2, g).
\]By definition of $\lambda_\infty$,  we know that
\[
\lambda(B(x_k, 2r, g),  g) >\lambda_\infty - \e
\]when $k$ is sufficiently large.  So we get
\[
\lambda(B(x_\infty,  r))  > \lambda_\infty - 2 \epsilon.
\] By (\ref{lamb>lambB}), we finally deduce
\[
\lambda =\lambda(\M) > \lambda_\infty - 3 \epsilon.
\] Since $\e$ can be sufficiently small, we have reached a contradiction with the assumption that
$\lambda<\lambda_\infty$. This shows that Case 2 can not happen, and only Case 1 occurs, implying that an extremal exists.

The bound for the extremal $v$ in the theorem is already proven in Lemma \ref{ledecay}.
This proves part (a) of the theorem.  \qed

\section{Proof of the theorem \ref{thmain} (b), the nonexistence part}

The proof is done by constructing a concrete 3 manifold on which the Log Sobolev functional does
not have an extremal. In order to present the main idea of the construction, we informally describe a crude example of a disconnected manifold of such kind. \\

{\it Example 3.1.} Let $(M_k, g_k)$, $k=1, 2, ...$, be a sequence of compact manifolds without
boundary and let $\lambda_k$ be the infimum of the Log Sobolev functional on $M_k$.
We assume that $\lambda_k$ is a strictly decreasing sequence bounded from below by a finite number. For instance we can take $M_k = (1+ k^{-2}) (S^1 \times S^1)$, the flat 2 torus whose
metric is the standard one scaled by the factor $1+ k^{-2}$.
Let $M$ be the disjoint union of $M_k$. We now prove that the Log Sobolev functional does
not have an extremal on $M$.
Suppose for contradiction that $v$ is an extremal of the Log Sobolev functional  on $M$, whose
infimum is $\lambda$. Then $\lambda<\lambda_k$ and
\be
\lab{L=L1+...}
\lambda = L(v, g) =\Sigma^\infty_{k=1} \int_{M_k} (4 |\nabla v|^2 +R_k v^2 - v^2 \ln v^2) dg_k.
\ee Here $R_k$ is the scalar curvature of $(M_k, g_k)$. Without loss of generality, we can assume
that $v|_{M_k}$ is not identically zero for $k=1, 2, 3, ...$.  Otherwise, we just delete those $M_k$
where $v|_{M_k}$ is identically zero.  Write
\[
v_k = \frac{v|_{M_k}}{\Vert v|_{M_k} \Vert_{L^2(M_k, g_k)}}.
\] Then, $\Vert v_k \Vert^2_{L^2(M_k, g_k)} = 1$ and
\[
\al
&\int_{M_k} (4 |\nabla v|^2 +R_k v^2 - v^2 \ln v^2) dg_k\\
&= \Vert v|_{M_k} \Vert^2 _{L^2(M_k, g_k)} \, \int_{M_k} (4 |\nabla v_k|^2 +R_k v_k^2 - v_k^2 \ln v_k^2) dg_k
   - \Vert v|_{M_k} \Vert^2_{L^2(M_k, g_k)} \ln  \Vert v|_{M_k} \Vert^2_{L^2(M_k, g_k)}\\
   &\ge  \Vert v|_{M_k} \Vert^2 _{L^2(M_k, g_k)} \, \int_{M_k} (4 |\nabla v_k|^2 +R_k v_k^2 - v^2_k \ln v^2_k) dg_k.
 \eal
 \]Here we used the fact that $ \Vert v|_{M_k} \Vert^2_{L^2(M_k, g_k)} \le \Vert v \Vert^2_{L^2(M)}
 =1$.  Hence
 \[
 \int_{M_k} (4 |\nabla v|^2 +R_k v^2 - v^2 \ln v^2) dg_k \ge \Vert v|_{M_k} \Vert^2 _{L^2(M_k, g_k)}
 \lambda_k.
 \]Substituting this to (\ref{L=L1+...}), we deduce
 \[
 \lambda \ge \Sigma^\infty_{k=1} \Vert v|_{M_k} \Vert^2 _{L^2(M_k, g_k)}
 \lambda_k.
 \]Notice that
 \[
 1 = \Vert v \Vert^2_{L^2(M)} = \Sigma^\infty_{k=1} \Vert v|_{M_k} \Vert^2 _{L^2(M_k, g_k)}.
 \]Multiplying this equality by $\lambda$ and subtracting the last inequality, we find that
 \[
  \Sigma^\infty_{k=1} \Vert v|_{M_k} \Vert^2 _{L^2(M_k, g_k)} (\lambda_k -\lambda)  \le 0,
 \] which is a contradiction with the fact that $\lambda_k >\lambda$.  Hence no such extremal
 $v$ exists.

 The manifold $M$ in this example is disconnected and therefore it can not serve as a proof of the
 theorem. However, building on the main idea from this example, we will construct
a manifold $\M$ which is a connected sum of infinitely many copies
of compact manifolds, each of which can be graphically described
as a ball with a handle or just a "hand bag". See the figure in
Step 4 of the proof.
 The basic components of the manifold are: round necks, truncated $S^3$, and tubes whose cross sections are the flat torus
 $S^1 \times S^1$.
By studying the behavior of the Log Sobolev functional when these
components are pasted together, we will eventually show that the
Log Sobolev functional does not have an extremal.

First let us introduce some notations.

\begin{definition}
\lab{decyltub}
(Round necks and flat tubes)

Let $h, A, B$ be real numbers, we use $N=N(h, A, B)$ to denote
the round neck $h^2 S^2 \times [A, B]$ with the  product metric $g=h^2 g_{S^2} \times g_{R^1}$.
Here $g_{S^2}$ is the standard round metric on $S^2$ with radius $1$;  $g_{R^1}$ is the Euclidean metric on $R^1$; and $h^2$ scales $g_{S^2}$ only. For convenience, we also normalize the scalar curvature
corresponding to $g_{S^2}$ to be $1$.
Let $x \in N(h, A, B)$. We use $x=(x_1, x_2, x_3)$ as a coordinate for $x$, where $(x_1, x_2) \in S^2$
and $x_3 \in [A, B]$.

If $A=0$, we will use $N(h, B)$ to denote $N(h, A, B)$.

We use $H=H(h, A, B)$ to denote the flat tube $h^2 (S^1 \times S^1) \times [A, B]$ with the product metric $g=h^2 g_{S^1 \times S^1} \times g_{R^1}$.
Here $g_{S^1 \times S^1}$ is the standard flat metric on $S^1 \times S^1$ so that the radius
of $S^1$ is $1$;  $g_{R^1}$ is the Euclidean metric on $R^1$; and $h^2$ scales $g_{S^1 \times S^1}$ only.
Let $x \in H(h, A, B)$. We use $x=(x_1, x_2, x_3)$ as a coordinate for $x$, where $(x_1, x_2) \in S^1 \times S^1$
and $x_3 \in [A, B]$.
\end{definition}

We need a number of lemmas again.

\begin{lemma}
\lab{lev<1/a3}
Let $v$ be a bounded, positive subsolution to the equation
(\ref{eqln}) in the round neck $N= h^2 S^2 \times [-l, l]$. i.e.
\[
4 \Delta v - R v + 2 v \ln v + \lambda v \ge 0.
\] Suppose $\lambda \le 0$,  $h \in (0, 1]$, $l \ge 2$ and
that $\Vert v \Vert_{L^2(N)} \le 1$. Then there exists a positive constant $C$ which is independent of $h$ such that
\[
 v^2(x) \le C \int_{B(x, 1)} v^2 dg
\] when $x \in  h^2 S^2 \times [-l+1, l-1]$.
\proof
\end{lemma}

The result in this lemma and the proof are analogous to that in
Lemma \ref{lemvi}. However, there is difference, namely the
constant $C$ in the lemma is independent of $h \in (0, 1]$.

First, we claim that there exists a positive constant $S_0$, independent of $h$, such that
 such that,
\be
\lab{sobcyl}
S_0 \left(\int u^{2n/(n-2)} dg \right)^{(n-2)/n} \le  \int (4 |\nabla u |^2 + R u^2)dg, \qquad n=3,
\ee for all $ u \in C^\infty_0(h^2 S^2 \times [-l,, l])$.
Here is a quick proof of the claim.
Consider the infinite round neck $S^2 \times h^{-2} R^1$.  Here $h^{-2} R^1$ is $R^1$equipped with
the scaled metric $h^{-2} g_{R^1}$. Note the curvature bounds and the lower bound of injectivity radius
are independent of $h$.  i.e. the necks have uniformly bounded geometry.  By \cite{Au:1}, there exists
a positive constant $S_0$ such that
\[
S_0 \left(\int u^{2n/(n-2)} dg \right)^{(n-2)/n} \le  \int (|\nabla u |^2 +  u^2)dg
\] for all $ u \in C^\infty_0(S^2 \times h^{-2} R^1)$. Notice that the scalar curvature of $S^2 \times h^{-2} R^1$ is
the constant $1$. Hence
\[
S_0 \left(\int u^{2n/(n-2)} dg \right)^{(n-2)/n} \le  \int (4 |\nabla u |^2 + R u^2)dg
\] for all $ u \in C^\infty_0(S^2 \times h^{-2} R^1)$. But this Sobolev inequality is scaling invariant.
Hence, for all $ u \in C^\infty_0( h^2S^2 \times R^1)$, inequality (\ref{sobcyl}) holds,
proving the claim.

Since $v$ is a subsolution of (\ref{eqln}) and $\lambda \le 0$ by assumption, given any $p \ge 1$, it is easy to see that
\[
- 4 \Delta v^p + p R v^p \le 2 p v^p \ln v.
\]
We select a smooth cut off function $\phi$ supported in $h^2 S^2 \times [-l,, l]$. Writing $w = v^p$ and using $w \phi^2$ as a test function in the above inequality, we deduce
\[
4 \int \nabla (w \phi^2) \nabla w dg + p \int R (w \phi)^2 dg \le 2 p \int (w \phi)^2 \ln v dg.
\] Since the scalar curvature is positive, this shows
\[
4 \int \nabla (w \phi^2) \nabla w dg + \int R (w \phi)^2 dg \le  2 p \int (w \phi)^2 \ln v^2 dg,
\]which induces, after integration by parts,
\[
 \int (4 |\nabla (w \phi)|^2 + R (w \phi)^2 )dg \le  4 \int |\nabla
\phi|^2 w^2 dg  + 2 p \int (w \phi)^2 \ln v^2 dg.
\]Applying  (\ref{sobcyl}) on the left hand side, we deduce
\[
S_0 \left(\int (w \phi)^{2n/(n-2)} dg \right)^{(n-2)/n} \le
 4 \int |\nabla
\phi|^2 w^2 dg  + 2 p \int (w \phi)^2 \ln v^2 dg.
\]

Now pick $x \in  h^2 S^2 \times [-l+1, l-1]$. Then $B(x, 1) \subset h^2 S^2 \times [-l+1, l-1]$.
Now we choose $\phi$ as suitable cut-off functions supported in $B(x, 1)$.
The rest of the proof of the lemma is the same as the proof of Lemma \ref{lemvi} after (\ref{9.5}), with
$\lambda$ there taken as $0$.
\qed

\medskip

The next lemma says that if $v$ is a solution of (\ref{eqln}) in a very long round neck, whose $L^2$ norm is less than $1$, then
$v$ is exponentially small in the middle section of the neck.

\begin{lemma}
\lab{lev<e-l}
There exists $h_0 \in (0, 1]$ such that the following statement holds for all $h \in (0, h_0]$.
Let $v$ be a smooth positive solution to the equation
(\ref{eqln}) in the round neck $N= h^2 S^2 \times [-l, l]$. Suppose $\lambda \le 0$, $l \ge 2$ and
that $\Vert v \Vert_{L^2(N)} \le 1$. Then there exist positive constants $a$ and $A$,  independent of $h$, such that
\[
  \int_{ h^2 S^2 \times [-l/2, l/2]} v^2 dg \le A e^{- a l} \, [  \int_{ h^2 S^2 \times [-l, -l+2]} v^2 dg
  +\int_{ h^2 S^2 \times [l-2, l]} v^2 dg ]
\] and
\[
v(x) \le  A e^{- a l}, \quad x \in  h^2 S^2 \times [-l/2, l/2].
\]
\proof
\end{lemma}

By the previous lemma, for $x  \in  h^2 S^2 \times [-l+1, l-1]$, we have a constant $C$ such that
\[
v(x) \le C.
\] Note the scalar curvature $R=1/h^2$. Hence there exists $h_0 \in (0, 1]$ such that if $h \in (0, h_0]$
then
\[
R/2 - 2 \ln v \ge 1/(2 h^2_0) - 2 \ln C \ge 0.
\] Combining this with equation (\ref{eqln}) i.e.
\[
4 \Delta v - R v + 2 v \ln v + \lambda v =0,
\] we find that  $v$ satisfies the inequality
\be
\lab{ddv-v}
\Delta v - \frac{1}{8 h^2_0} v  \ge 0 \quad \text{in} \quad h^2 S^2 \times [-l+1, l-1].
\ee Here we have used the assumption that $\lambda \le 0$.

We pick a cut off function $\phi \in C^\infty_0(N)$, satisfying the following requirements.
\[
\al
\phi(x)=\phi(x_1, x_2, x_3)=
\begin{cases}
0, &\quad x_3 \in [-l, -l+1] \cup [l-1, l],\\
\text{a number in } \quad (0, 1), &\quad x_3 \in [-l+1, -l+2] \cup [l-2, l-1]\\
1, &\quad x_3 \in [-l+2, l-2].\\
\end{cases}
\eal
\] We also require that $|\nabla \phi| \le 4$. Here we recall that $x_3$ is the longitudinal component of
the coordinate of the point $x$ in the neck $N$, as described in
Definition \ref{decyltub}. See the figure below. \\

 \psfrag{X}{$X$}
 \psfrag{p0}{$\phi=0$}
 \psfrag{p01}{$\phi \in (0, 1)$}
 \psfrag{p1}{$\phi=1$}
 \psfrag{N}{$N$}
 \psfrag{l}{$l$}
 \psfrag{-l}{$-l$}
 \psfrag{-l1}{$-l+1$}
 \psfrag{l-1}{$l-1$}
 \centerline{\epsfig{figure=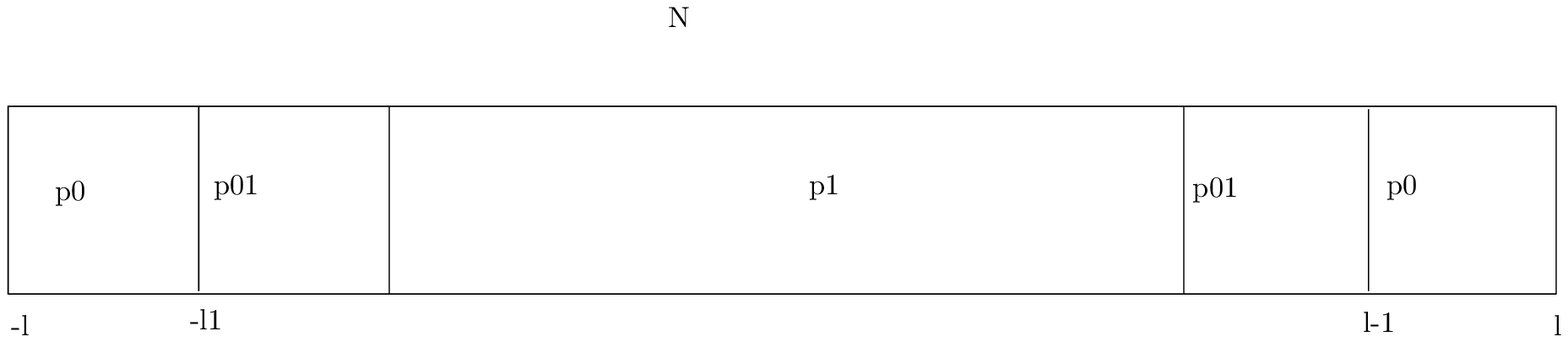, height=2cm,
 width=14cm}}

Let $a$ be a positive number to be determined later. Using $e^{a (l-|x_3|)} \phi^2 v$ as a test function on
(\ref{ddv-v}) and performing integration by parts, we find that
\be
\lab{iy123}
\al
& \frac{1}{8 h^2_0} \int e^{a (l-|x_3|)} \phi^2 v^2 dg \le  \int e^{a (l-|x_3|)} \phi^2 v \Delta v dg \\
&= - \int e^{a (l-|x_3|)} \phi^2 |\nabla v|^2 dg - 2  \int e^{a (l-|x_3|)} v \phi \nabla \phi \nabla v dg
- \int e^{a (l-|x_3|)}  \nabla (a (l-|x_3|))  \nabla v v \phi^2 dg\\
&\equiv -Y_1 - Y_2 -Y_3.
\eal
\ee We need to bound $|Y_2|$ and $|Y_3|$.

First we notice
\[
\al
|Y_2| &\le 2 \int e^{a (l-|x_3|)} v \phi |\nabla \phi \nabla v| dg\\
&\le \frac{1}{4}  \int e^{a (l-|x_3|)} \phi^2 |\nabla v|^2 dg +4  \int e^{a (l-|x_3|)} v^2 |\nabla \phi |^2 dg.
\eal
\]Therefore
\be
\lab{y2<}
|Y_2| \le  \frac{1}{4} Y_1 + 4 \int_{supp \nabla \phi}  e^{a (l-|x_3|)} v^2 dg.
\ee

Next
\[
\al
|Y_3| &\le a \int e^{a (l-|x_3|)}  |\nabla v|  v \phi^2 dg\\
&\le \frac{a}{2} \int e^{a (l-|x_3|)}   \phi^2 v^2 dg + \frac{a}{2} \int e^{a (l-|x_3|)}  |\nabla v|^2  \phi^2 dg\\
&=\frac{a}{2} \int e^{a (l-|x_3|)}   \phi^2 v^2 dg + \frac{a}{2} Y_1.
\eal
\]Choosing $a \le 1$ and substituting this and (\ref{y2<}) into (\ref{iy123}), we deduce
\[
\frac{1}{8 h^2_0} \int e^{a (l-|x_3|)} \phi^2 v^2 dg \le
4 \int_{supp \nabla \phi}  e^{a (l-|x_3|)} v^2 dg + \frac{a}{2} \int e^{a (l-|x_3|)}   \phi^2 v^2 dg.
\]Taking $a = \min \{1, \frac{1}{8 h^2_0} \}$, we arrive at
\be
\lab{eav2<}
\int e^{a (l-|x_3|)} \phi^2 v^2 dg \le C h^2_0 \int_{supp \nabla \phi}  e^{a (l-|x_3|)} v^2 dg.
\ee Observe that  when $ x \in supp \nabla \phi$ we have
\[
0 \le l - |x_3| \le 2.
\] Also, when $x \in h^2 S^2 \times [-2l/3, 2l/3]$, i.e. when $-2l/3 \le x_3 \le 2l/3$, we have
\[
l - |x_3| \ge l/3, \qquad \phi(x)=1.
\]Therefore (\ref{eav2<}) implies
\[
 \int_{ h^2 S^2 \times [-2l/3, 2l/3]} v^2 dg \le C h^2_0 e^{2 a} e^{- a l/3} \, [  \int_{ h^2 S^2 \times [-l, -l+2]} v^2 dg
  +\int_{ h^2 S^2 \times [l-2, l]} v^2 dg ]
\]which yields the desired integral bound, after adjusting the coefficients.  The pointwise  bound in the lemma
is an immediate consequence the integral bound and Lemma \ref{lev<1/a3}
\qed  \\

 Let $v$ again be a positive solution of (\ref{eqln}) in a very long round neck, whose $L^2$ norm is less than $1$. The next lemma says that if $v$ vanishes at one end of the neck, then
$v$ is exponentially small near that end.

\begin{lemma}
\lab{2lev<e-l}
There exists $h_0 \in (0, 1]$ such that the following statement holds for all $h \in (0, h_0]$.
Let $v$ be a smooth positive solution to the equation
(\ref{eqln}) in the round neck $N= h^2 S^2 \times [0, l]$. Suppose $\lambda \le 0$, $l \ge 2$ and
that $\Vert v \Vert_{L^2(N)} \le 1$.  Suppose also $v(x)=0$ when $x \in h^2 S^2 \times \{ l \}$. i.e.
$v$ vanishes at the right end of the neck. Then there exist positive constants $a$ and $A$,  independent of $h$, such that
\[
  \int_{ h^2 S^2 \times [l/2, l]} v^2 dg \le A e^{- a l} \int_{h^2 S^2 \times [0, 1]} v^2 dg.
\]
\proof
\end{lemma}

We extend $v=v(x)$ to a function on the longer neck $h^2 S^2 \times [0, l+1]$ by
assigning $v(x)=0$ when $x_3 \ge l$. Since $v(x)=0$ when $x_3 =l$, it is easy to see that the
extended $v$ is a subsolution to (\ref{eqln}) on $h^2 S^2 \times [0,  l+1]$.
By Lemma \ref{lev<1/a3}, for $x  \in  h^2 S^2 \times [1, l]$, there exists a constant $C$ such that
\[
v(x) \le C.
\] Since the scalar curvature $R=1/h^2$,  there exists $h_0 \in (0, 1]$ such that if $h \in (0, h_0]$
then
\[
R/2 - 2 \ln v \ge 1/(2 h^2_0) - 2 \ln C \ge 0.
\] Combining this with equation (\ref{eqln}) i.e.
\[
4 \Delta v - R v + 2 v \ln v + \lambda v =0.
\] we find that  $v$ satisfies the inequality
\be
\lab{ddv-vl}
\Delta v - \frac{1}{8 h^2_0} v  \ge 0 \quad \text{in} \quad h^2 S^2 \times [1, l].
\ee Here we have again used the assumption that $\lambda \le 0$.

We pick a cut off function $\phi \in C^\infty_0(N)$, satisfying $|\nabla \phi| \le 4$ and the following requirements.
\[
\al
\phi(x)=\phi(x_1, x_2, x_3)=
\begin{cases}
0, \quad &x_3 \in [0, 1],\\
\text{a number in } \quad (0, 1), \quad &x_3 \in [1, 2] \\
1, \quad &x_3 \in [2, l].
\end{cases}
\eal
\]

Let $a$ be a positive number to be determined later. Using $e^{a x_3} \phi^2 v$ as a test function on
(\ref{ddv-vl}) and performing integration by parts, we find that
\be
\lab{2iy123}
\al
& \frac{1}{8 h^2_0} \int e^{a x_3} \phi^2 v^2 dg \le  \int e^{a x_3} \phi^2 v \Delta v dg \\
&= - \int e^{a x_3} \phi^2 |\nabla v|^2 dg - 2  \int e^{a x_3} v \phi \nabla \phi \nabla v dg
- \int e^{a x_3}  \nabla (a x_3)  \nabla v v \phi^2 dg\\
&\equiv -Y_1 - Y_2 -Y_3.
\eal
\ee Note that boundary terms vanish since $v=0$ at the right end of the neck and
$\phi=0$ at the left end.  Let us bound $|Y_2|$ and $|Y_3|$.

First we notice
\[
\al
|Y_2| &\le 2 \int e^{a x_3} v \phi |\nabla \phi \nabla v| dg\\
&\le \frac{1}{4}  \int e^{a x_3} \phi^2 |\nabla v|^2 dg + 4 \int e^{a x_3} v^2 |\nabla \phi |^2 dg.
\eal
\]Therefore
\be
\lab{2y2<}
|Y_2| \le  \frac{1}{4} Y_1 + 4 \int_{supp \nabla \phi}  e^{a x_3} v^2 dg.
\ee

Next
\[
\al
|Y_3| &\le a \int e^{a x_3}  |\nabla v|  v \phi^2 dg\\
&\le \frac{a}{2} \int e^{a x_3}   \phi^2 v^2 dg + \frac{a}{2} \int e^{a x_3}  |\nabla v|^2  \phi^2 dg\\
&=\frac{a}{2} \int e^{a x_3}   \phi^2 v^2 dg + \frac{a}{2} Y_1.
\eal
\]Choosing $a \le 1$ and substituting this and (\ref{2y2<}) into (\ref{2iy123}), we deduce
\[
\frac{1}{8 h^2_0} \int e^{a x_3} \phi^2 v^2 dg \le
4 \int_{supp \nabla \phi}  e^{a x_3} v^2 dg + \frac{a}{2} \int e^{a x_3}   \phi^2 v^2 dg.
\]Taking $a = \min \{1, \frac{1}{8 h^2_0} \}$, we arrive at
\be
\lab{2eav2<}
\int e^{a x_3} \phi^2 v^2 dg \le C h^2_0 \int_{supp \nabla \phi}  e^{a x_3} v^2 dg.
\ee Observe that  when $ x \in supp \nabla \phi$ we have
\[
0 \le x_3 \le 1.
\] Also, when $x \in h^2 S^2 \times [l/2, l]$, we have
\[
x_3 \ge l/2, \qquad \phi(x)=1.
\]Therefore (\ref{2eav2<}) implies
\[
 \int_{ h^2 S^2 \times [l/2, l]} v^2 dg \le C h^2_0 e^a e^{- a l/2}\int_{h^2 S^2 \times [0, 1]} v^2 dg,
\]proving the lemma.
\qed

\medskip

The following lemma is similar to Lemma \ref{leLMD}. The difference is that we are comparing the
infimum of the Log Sobolev functionals on two different domains in this lemma. The proof is almost
identical.

\begin{lemma}
\label{leLXY} Let $E$ and $F$ be two  domains of $\M$ such that $E
\subset F$ and that $E$ is compact. Let $v \in W^{1, 2}_0(F)$,
$\Vert v \Vert_{L^2(F)}=1$ be an extremal of $\lambda(F)$ so that
it is a smooth positive solution  of the equation
 \[
 4 \Delta v - R v + 2 v \ln v + \lambda(F) v=0.
 \]  For any smooth cut-off
function $\eta$ such that $\eta v \in C^\infty_0(E)$ and $0 \le
\eta \le 1$, it holds
\[
\lambda (E) \le \lambda(F) + 4 \frac{ \int  v^2 |\nabla \eta|^2
dg}{\int (v \eta)^2 dg} - \frac{\int (v \eta)^2 \ln \eta^2
dg}{\int (v \eta)^2 dg}.
\]
\proof
\end{lemma}

Since $\eta v/\Vert \eta v \Vert_2 \in C^\infty_0(E)$  and its
$L^2$ norm is $1$, we have, by  definition,
\[
\lambda(E) \le \int \left[ 4 \frac{|\nabla (\eta v)|^2 }{\Vert
\eta v \Vert^2_2 } + R \frac{(\eta v)^2}{\Vert \eta v \Vert^2_2 }
- \frac{(\eta v)^2}{\Vert \eta v \Vert^2_2 } \ln \frac{(\eta
v)^2}{\Vert \eta v \Vert^2_2 } \right] dg.
\]This implies
\begin{equation}
\lambda(E) \Vert \eta v \Vert^2_2  \le \int \left[ 4 |\nabla (\eta
v)|^2 +R (\eta v)^2 -
 (\eta v)^2 \ln (\eta v)^2 \right] dg  + \Vert \eta v \Vert^2_2 \ln  \Vert \eta v \Vert^2_2.
 \label{3.14}
\end{equation}

 On the other hand, $v$ is a smooth positive solution  of the equation
 \[
 4 \Delta v - R v + 2 v \ln v + \lambda(F) v=0.
 \]Using $\eta^2 v$ as a test function for the equation, we find
 \[
 \lambda (F) \int (\eta v)^2 dg = - 4 \int (\Delta v) \eta^2 v dg + \int R (\eta v)^2 dg -
 2 \int (\eta v)^2 \ln v dg.
 \]Using integration by parts, we deduce
 \[
 -4 \int (\Delta v)  \eta^2 v dg= 4 \int |\nabla (\eta v)|^ 2 dg- 4 \int v^2 |\nabla \eta|^2 dg.
 \]Hence
 \begin{equation}
  \lambda(F) \int (\eta v)^2 dg= 4 \int |\nabla (\eta v)|^ 2 dg- 4 \int v^2 |\nabla \eta|^2 dg
 + \int R (\eta v)^2 dg-
 2 \int (\eta v)^2 \ln v dg.
 \label{3.15}
 \end{equation} Comparing (\ref{3.15}) with (\ref{3.14}) and noting that $\Vert \eta v \Vert_2 \le 1$, we
 obtain
 \[
 \lambda(E) \Vert \eta v \Vert^2_2  \le  \lambda(F) \Vert \eta v \Vert^2_2 +
 4 \int |\nabla \eta|^2 v^2 dg - \int (\eta v)^2 \ln \eta^2 dg.
 \]
 \qed \\

 The following lemma says that if a domain $E$ contains a round neck of length $l$ and
 $F$ is the extension of $E$, which is obtained by pasting a segment of the round neck with length $1$,
 then $|\lambda(E)-\lambda(F)|$ is exponentially small. \\

\begin{lemma}
\lab{leLX-LY} Let $E \subset \M$ be a compact domain such that
\[
E =X_0 \cup N(h, l)
\]which is the connected, non-overlapping  union of a domain $X_0$ with the round neck $N(h, l) = h^2 S^2 \times [0, l]$.
Let
\[
F =X_0 \cup N(h, l+1)
\]which is the connected, non-overlapping union of
 $X_0$ with the round neck $N(h, l+1) = h^2 S^2 \times [0, l+1]$.
There is $h_0 \in [0, 1]$ and $l_0>0$
such that for all $h \in [0, h_0]$ and $l \ge l_0$,  the following statement holds:

If $\lambda(E) \le 0$, then there exist positive numbers $a$ and
$A$ such that
\[
\lambda(F) \ge \lambda(E) - A e^{-a l}.
\]
\proof
\end{lemma}

First let us see the figure depicting $E$ and $F$ below.

\psfrag{X0}{$X_0$} \psfrag{Nhl}{$N(h, l)$}
\psfrag{Nhl1}{$N(h,l+1)$} \psfrag{F}{$F=X_0 \cup N(h, l+1)$}
\psfrag{E}{$E=X_0 \cup N(h, l)$}
 \centerline{\epsfig{figure=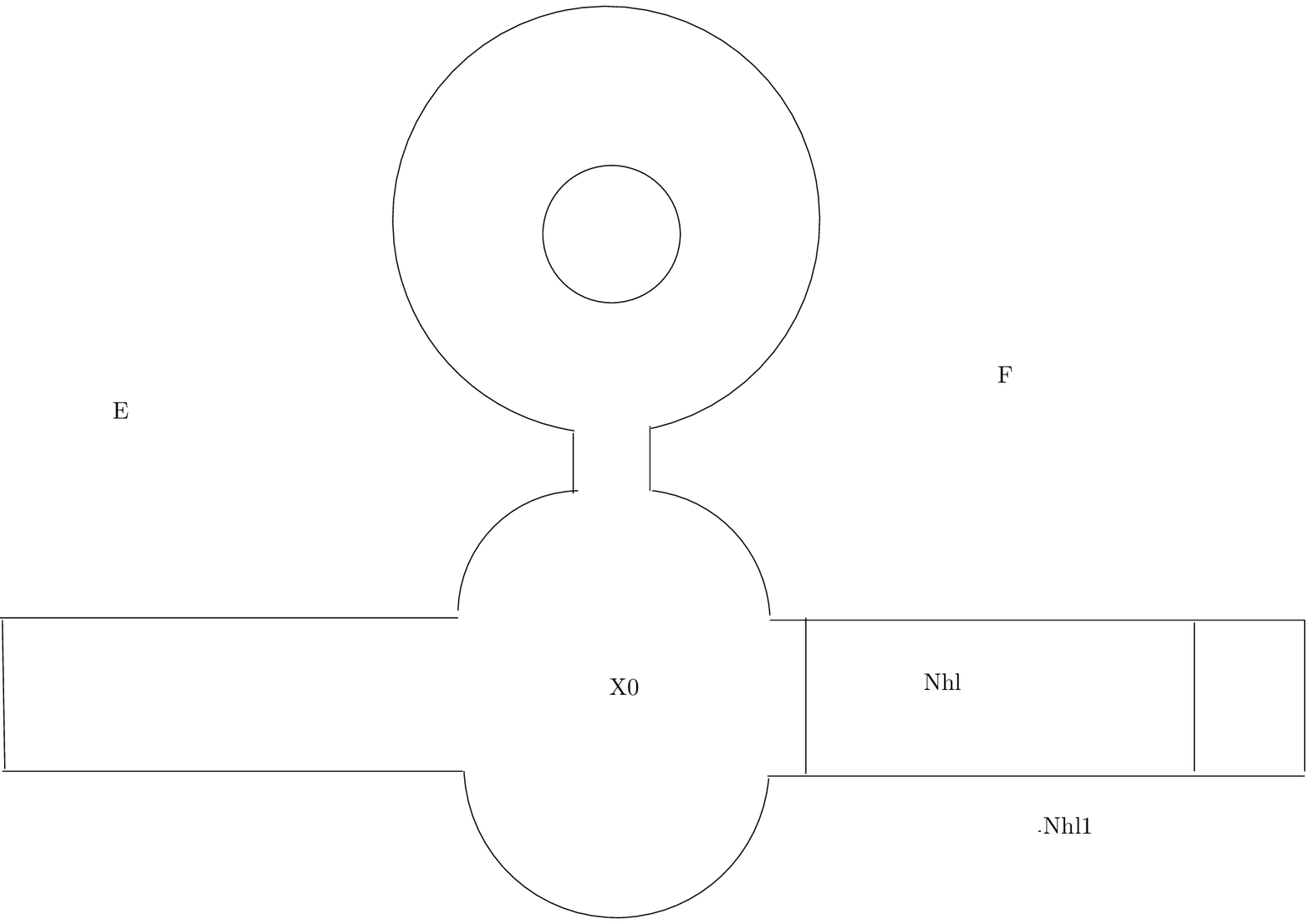, height=5cm, width=14cm}}

Pick a smooth cut off function $\eta$ such that $|\nabla \eta| \le 4$ and that
\[
\al
\eta =\eta(x)=
\begin{cases}
0, \qquad &x \in h^2 S^2 \times [l, l+1]\\
\text{a number in} \quad (0, 1), \qquad &x \in h^2 S^2 \times [l-1, l]\\
1, \qquad &x \in F -  (h^2 S^2 \times [l-1, l+1]).
\end{cases}
\eal
\] Let $v$ be an extremal for $\lambda(F)$, which exists since $F$ is compact. Then
$\eta v \in C^\infty_0(E)$.  By Lemma \ref{leLXY}, we have
\be
\lab{laxlay} \lambda (E) \le \lambda(F) + 4 \frac{ \int  v^2
|\nabla \eta|^2 dg}{\int (v \eta)^2 dg} - \frac{\int (v \eta)^2
\ln \eta^2 dg}{\int (v \eta)^2 dg}. \ee
Observe that
\[
\int (v \eta)^2 dg = \int v^2 dg - \int v^2 (1-\eta^2) dg \ge 1 -
\int_{h^2 S^2 \times [l-1, l+1]} v^2 dg.
\]Using Lemma \ref{2lev<e-l} on $h^2 S^2 \times [0, l+1]$, we infer, for some positive numbers $a$ and $A$, that
\[
\int_{h^2 S^2 \times [l-1, l+1]} v^2 dg \le A e^{-a l}.
\]Hence
\[
\int (v \eta)^2 dg \ge 1- A e^{-a l}.
\]Also notice that
\[
\int  v^2 |\nabla \eta|^2 dg \le 16 \int_{h^2 S^2 \times [l-1, l]}  v^2  dg \le 16 A e^{-a l},
\]and
\[
| \int (v \eta)^2 \ln \eta^2 dg| \le e^{-1}  \int_{h^2 S^2 \times [l-1, l]}  v^2  dg \le A e^{-a l}.
\] Substituting the last three inequalities into (\ref{laxlay}), we deduce
\[
\lambda (E) \le \lambda(F) + C \frac{ A e^{-a l}}{1-A e^{-a l}}.
\]Therefore, there exists $l_0>0$ such that for all $l \ge l_0$, we have
\[
\lambda (E) \le \lambda(F) +A e^{-a l}
\]for some constant $A>0$, whose value may have been adjusted from the last line. \qed \\

The following lemma says that the infimum of the Log Sobolev functional on
a flat tube goes to $-\infty$ when the cross section of the tube goes to $0$.\\

\begin{lemma}
\lab{leLftube}
Let $H=H(h, 0, 1)=h^2(S^1 \times S^1) \times [0, 1]$ be a flat tube given in Definition \ref{decyltub}.
Then $\lambda(H(h, 0, 1)) \to -\infty$ when $h \to 0$.
\proof
\end{lemma}

Given $x \in H(h, 0, 1)$, let $(x_1, x_2, x_3)$ be its coordinate described in Definition \ref{decyltub}.
Consider the one variable function
\[
\al
v=v(x_3)=
\begin{cases}
4 \frac{\sqrt{3}}{\sqrt{8} \,  \pi h} x_3, &\qquad x_3 \in [0, 1/4],\\
\frac{\sqrt{3}}{\sqrt{8}  \,  \pi h} , &\qquad x_3 \in [1/4, 3/4],\\
\frac{\sqrt{3}}{\sqrt{8}  \,  \pi h} [1- 4(x_3 - 3/4)], & \qquad x_3 \in [3/4, 1].\\
\end{cases}
\eal
\] We compute
\[
\int_{H(h, 0, 1)} v^2 dg = 4 \pi^2 h^2 \int^1_0 v^2 dx_3
=\frac{3}{8 \pi^2 h^2} 4 \pi^2 h^2 ( 2 \int^{1/4}_0 16 x^2_3 dx_3 + \frac{1}{2} ) =1,
\]
\[
\int_{H(h, 0, 1)} |\nabla v|^2 dg = 4 \pi^2 h^2 \int^1_0 |\partial_{x_3} v|^2 dx_3
=\frac{3}{8 \pi^2 h^2} 4 \pi^2 h^2 ( 2 \int^{1/4}_0 16 dx_3 ) =12,
\]
\[
\al
\int_{H(h, 0, 1)} &v^2 \ln v^2 dg = 4 \pi^2 h^2 \int^1_0  v^2 \ln v^2 dx_3\\
&=\frac{3}{8 \pi^2 h^2} 4 \pi^2 h^2 [ 2 \int^{1/4}_0 16 x^2_3 \ln ( \frac{3}{8 \pi^2 h^2} 16 x^2_3) dx_3
+ \int^{3/4}_{1/4} \ln ( \frac{3}{8 \pi^2 h^2} ) dx_3] \\
&= - \frac{3}{4} \ln h^2 +c
\eal
\]where $c$ is a constant independent of $h$.

Since the scalar curvature is zero, these computation imply
\[
\lambda(H(h, 0, 1)) \le \int_{H(h, 0, 1)} (4 |\nabla v|^2 - v^2 \ln v^2) dg =  \frac{3}{4} \ln h^2 + c.
\]This shows $\lambda(H(h, 0, 1)) \to -\infty$ when $h \to 0$. \qed

Now we are ready to give
\medskip

{\bf Proof of Theorem \ref{thmain} (b).}

As mentioned earlier we will construct a noncompact manifold with bounded geometry such that
the Log Sobolev functional does not have an extremal.  The manifold is a connected sum
of infinitely many components connected by increasingly long round necks.
Each of the component shapes like a hand bag.  The handle of a hand bag is a flat tube of
certain  thickness.  By pinching the handle, we can control precisely the difference between the infimums of
the Log Sobolev functional on two adjacent hand bags.  The long round necks serve the following purpose: when two hand bags are joined, the change in the infimum
of the Log Sobolev functional happens in a controlled way.
 In the next few steps we will construct the components inductively in detail.
\medskip

{\it Step 1.} constructing the central component $\Omega_0$.
See the figure at the end of the step. \\

{\it Step 1.1.}  We start with the standard $3$ sphere with three small balls cut out.
To be more precise, let
\[
D=S^3-(B_1 \cup B_2 \cup B_3)
\] where $S^3$ is the standard $3$ sphere and $B_i=B(m_i, r)$, $i=1, 2, 3$, are geodesic balls on $S^3$
with radius $r>0$.  We take $m_1$, the center of the ball $B_1$ at
the north pole of $S^3$; $x_2$, the center of the ball $B_2$ at
the "left end" of the equator; and $x_3$, the center of the ball
$B_3$ at the "right end" of the equator. The radius $r$ is so
chosen that $\partial B_i$, $i=1, 2, 3$, is  $h^2 S^2$, the
standard $2$ sphere with radius $h$. The radius $h \in (0, 1/4]$
is made sufficiently small so that the following conditions hold:

(1)  Lemmas \ref{lev<e-l} and \ref{leLX-LY} hold;

(2)  $\lambda( h^2  (S^1 \times S^1) \times [-2, 2])) <0$. That is the infimum of the
Log Sobolev functional for the flat tube is negative.

By  Lemma \ref{leLftube}, condition (2) can always be satisfied when $h$ is small enough.

Once chosen, this $h$ will be fixed through out the proof. \\

{\it Step 1.2.} Attach a long round neck $h^2 S^2 \times [0,
l]$ to $D$ along $\partial B_2$ and $\partial B_3$ respectively.
Here $l>0$ is a large number given by \be \lab{l=} l=\max \{ l_0,
\frac{1}{a} \ln (1000 A e^{2a}/a^2),  \frac{1}{a} \ln (1000 e^{2 a} A), 2 \}.
\ee Here $l_0, a, A$ are the numbers in Lemmas \ref{lev<e-l},
\ref{2lev<e-l} and \ref{leLX-LY}. By taking this value for $l$,
all these three lemmas hold and
\be
 \lab{Ae-a<1/k2}
 20 A e^{2 a}  e^{-a(l+k)} \le \frac{1}{2(1+k^2)}, \quad k=0, 1, 2, 3, .... \ee  This inequality, to be used shortly in the end of the proof, can be
verified easily by finding the maximum of $(1+k^2) e^{-a k}$.

Let $h^2 (S^1 \times S^1) \times S^1 = h^2 (S^1 \times S^1) \times [-\pi, \pi]$ be a
flat $3$ torus, which is regarded as a flat tube given in Definition \ref{decyltub}.
Consider
\[
E = h^2 (S^1 \times S^1) \times [-\pi, \pi]-B_4.
\]Here $B_4=B(m_4, h)$
is the geodesic ball of radius $h$ centered at
$m_4$ whose coordinate is  $(0, 0, \pi)$. i.e. $m_4$ is at the bottom of the flat tube.
Note $h$ is less than the injectivity radius of the flat torus, which is $\pi h$.
Therefore we know $B_4$ is
isometric to the Euclidean ball of radius $h$. Hence $\partial B_4 = h^2 S^2$.

Now we join $D$ with $E$ by a short round neck $h^2 S^2 \times [0, 1]$
by pasting $h^2 S^2 \times \{0\}$ with $\partial B_1$, and pasting $h^2 S^2 \times \{0\}$ with $\partial B_4$. \\

{\it Step 1.3.} The metric near the pasted boundaries are smoothed out to satisfy the following conditions.

(1)  only the original metric on $D$ near a small neighborhood of $\partial B_i$, $i=1, 2, 3$, are perturbed, so that the metric on the attached long round necks stay the same.

(2) only the metric in a small neighborhood of $\partial B_4$ is perturbed so that
the metric on  $ h^2 (S^1 \times S^1) \times [-2, 2]$, which is the top portion of the flat tube,
 stays intact.

 Note the smoothing process is a standard procedure in geometry when one constructs connected
 sums of two manifolds.

The resulting manifold with boundary is called $\Omega_0$ with
metric $g_0$. By condition (2) in Step 1.1, we have
 \be
\lab{LOm1<0} \lambda(\Omega_0, g_0) \le \lambda( h^2  (S^1
\times S^1) \times [-2, 2])) <0. \ee

For clarity, we write \be \lab{om1=z1} \Omega_0 = Z_0 \cup X \cup
H \cup Y_0. \ee Here $Z_0$ is the round neck at the left,
which is $h^2 S^2 \times [0, l]$; $Y_0$ is the round neck at
the right, which is $h^2 S^2 \times [0, l]$ again. In order to
distinguish the two, we use $z$ to denote points in $Z_0$, and use
$y$ to denote points in $Y_0$.  $H$ denotes the top portion of the
flat tube where the third variable of the coordinates is in the
interval $[-2, 2]$. i.e. $h^2 (S^1 \times S^1) \times [-2, 2]$. We
will use the following global coordinate to denote the topological
$H$ in the rest of the proof.
\be \lab{H=} H=[-\pi, \pi]^2 \times
[-2, 2].
\ee The metric $g_0$ on $H$ is just $h^2 g_{S^1 \times
S^1} \times g_{R^1}$. The region $X$ is defined to be
\[
X =\Omega_0 - (Z_0 \cup H \cup Y_0)
\] We call $X$ the core of $\Omega_0$. The manifold $(X, g_0)$ will
serve as the core for all the rest of the components $\Omega_k$. \\

The shape of $\Omega_0$ is illustrated here.


 \psfrag{X}{$X$}
 \psfrag{Z}{$Z_0$}
 \psfrag{Y}{$Y_0$}
 \psfrag{H}{$(H, g_0)$}
 \psfrag{O0}{$\Omega_0$}
 \psfrag{zl}{length $l$}
 \psfrag{yl}{length $l$}
\centerline{\epsfig{figure=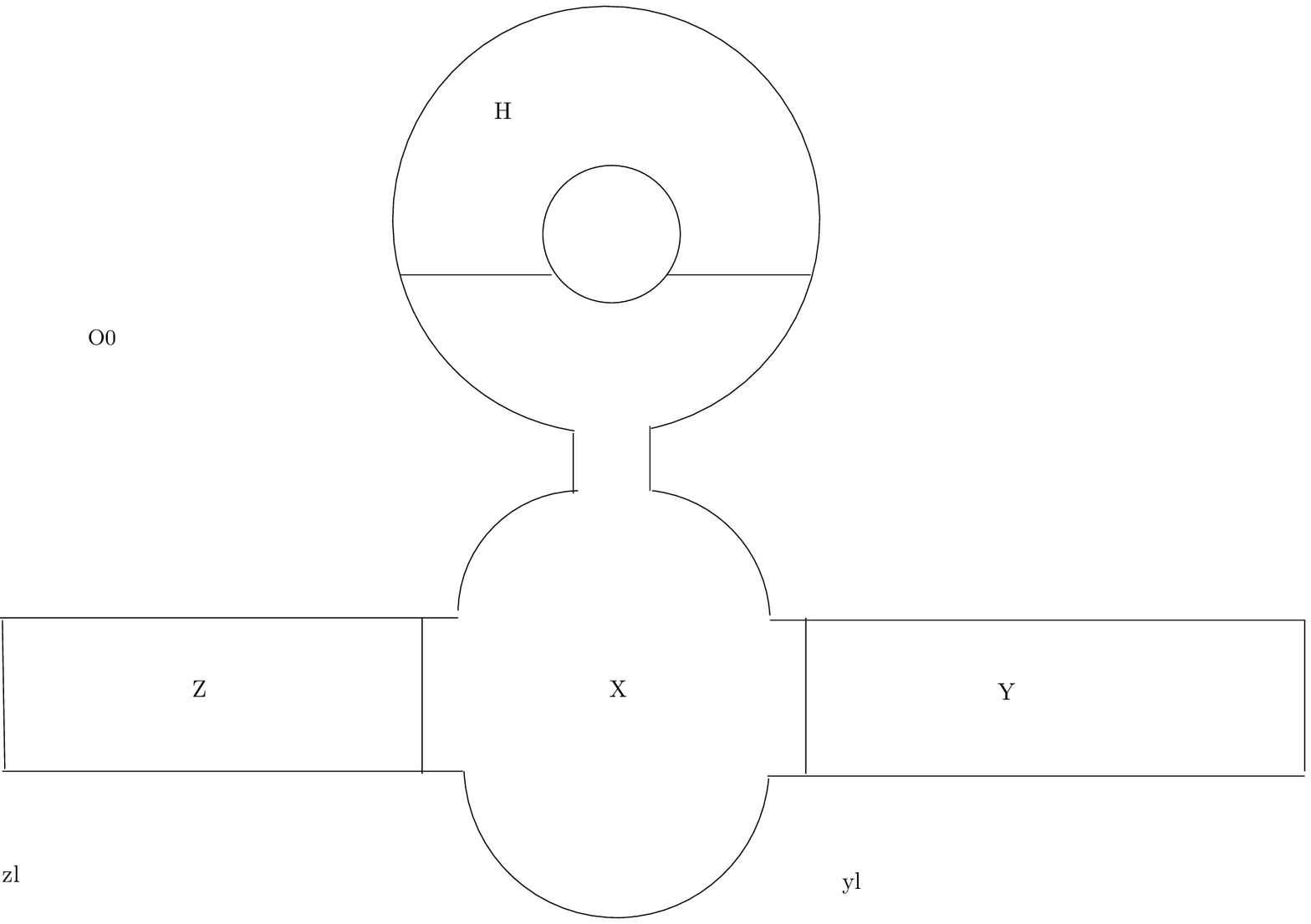, height=5cm}}

{\it Step 2.} constructing the next component $\Omega_1$ with
metric $g_1$ such that \be \lab{LOm2} \lambda(\Omega_1, g_1) =
\lambda(\Omega_0, g_0) -1. \ee \\

{\it Step 2.1.} Attach the round neck $h^2 S^2 \times [0, 1]$ to
the left end of $\Omega_0$, forming the round neck $h^2 S^2 \times
[0, l+1]$ on the left side, which we call $Z_1$. Then attach the
round neck $h^2 S^2 \times [0, 1]$ to the right end of $\Omega_0$,
forming the round neck $h^2 S^2 \times [0, l+1]$ on the right
side, which we call $Y_1$. The resulting domain is called
$\Omega_1$ with inherited metric called $\tilde g_1$. For
convenience we write
\[
\Omega_1 =  Z_1 \cup X \cup H \cup Y_1.
\] Note $\tilde g_1$ is already a smooth metric. In fact $\tilde g_1$ is the same as $g_0$ on $X$ and $H$,
and it is just the product metric on $h^2 g_{S^2} \times g_{R^1}$ on $Z_1$ and $Y_1$. But it is not
the desired one yet. \\

{\it Step 2.2.} Modify $\tilde g_1$ to a new metric $g_1$ so that
(\ref{LOm2}) holds. This modification only happens on $H$, the top
portion of the flat tube. More precisely, this is done by pinching
the top portion of the flat tube. Here are the details.

Recall that the top portion of $\Omega_1$ is the flat tube $H=
[-\pi, \pi]^2 \times [-2, 2]$. Let $\theta$ be a smooth function
on $\Omega_1$, satisfying
\[
\al
\theta=\theta(x)=
\begin{cases}
1, \qquad &x \in \Omega_1 - H\\
\text{a number in} \quad (1/2, 1), \qquad &x \in H, \quad x \in [-\pi, \pi]^2 \times [-2, -1]\\
1/2, \qquad &x \in H, \quad x \in [-\pi, \pi]^2 \times [-1, 1]\\
\text{a number in} \quad (1/2, 1), \qquad &x \in H, \quad x \in
[-\pi, \pi]^2 \times [1, 2].
\end{cases}
\eal
\] See the figure below.\\

\psfrag{th1}{$\theta=1$}
 \psfrag{th12}{$1/2 \le \theta \le 1$}
 \psfrag{th0.5}{$\theta=1/2$}
 \psfrag{H}{$H$}
 \psfrag{-2}{$x_3=-2$}
 \psfrag{-1}{$x_3=-1$}
 \psfrag{1}{$x_3=1$}
 \psfrag{2}{$x_3=2$}
\centerline{\epsfig{figure=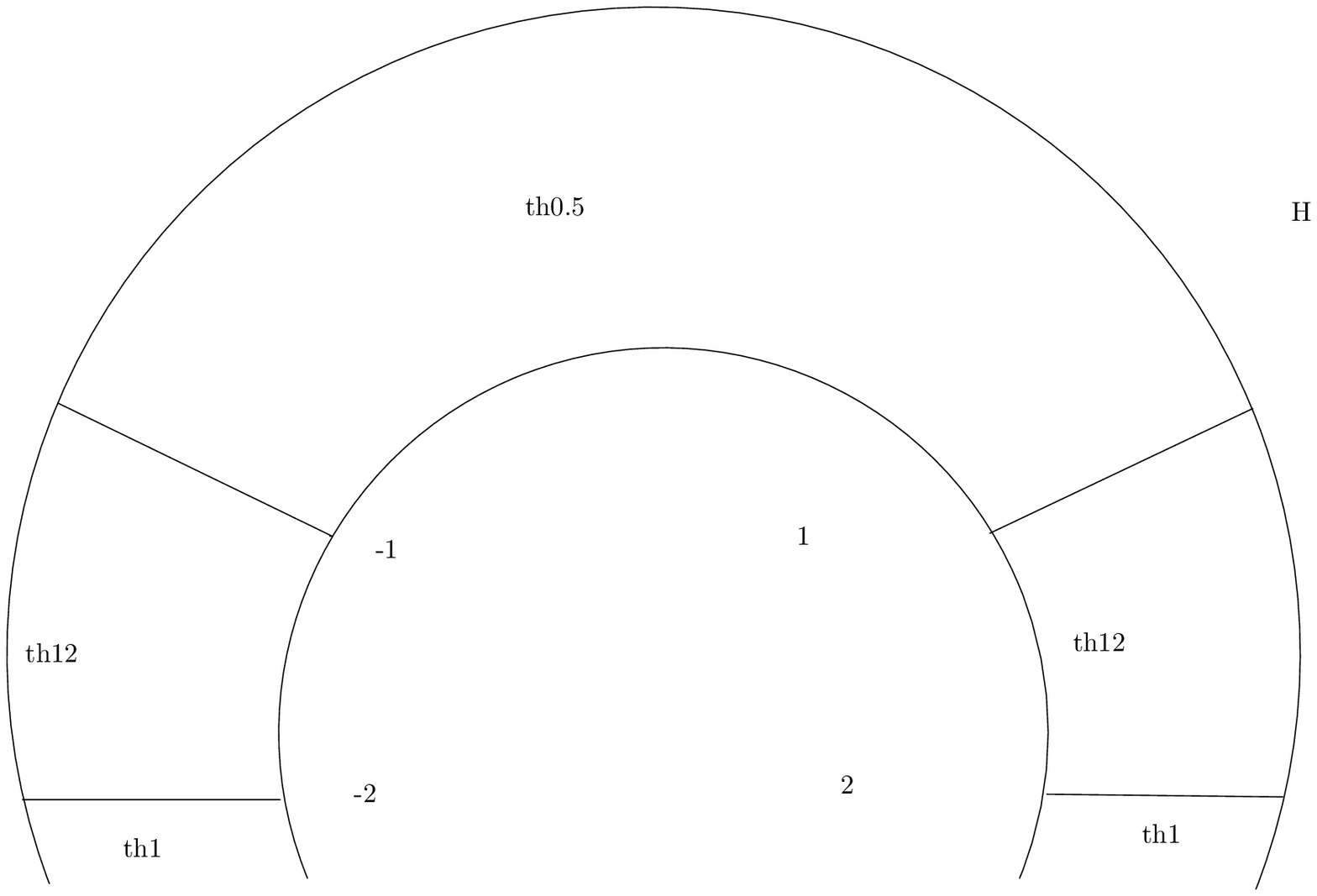, height=4cm, width=12cm}}

Now consider the metrics on $\Omega_1$:
\[
\al g^{(p)}_1(x)=
\begin{cases}
\tilde g_1(x), \qquad & x \in \Omega_1 - H\\
[\theta^p(x) h^2 g_{S^1 \times S^1}] \times g_{R^1},  \qquad
& x \in H.
\end{cases}
\eal
\]We claim that there exists a number $p_1>0$  so that
\be \lab{L=L-1} \lambda(\Omega_1, g^{(p_1)}_1) = \lambda(\Omega_0,
g_0) -1. \ee  Here is the proof. Regarding $(\Omega_0, g_0)$ as a
domain in $(\Omega_1, \tilde g_1)$ and applying Lemma
\ref{leLX-LY} twice, we know that
\[
\lambda(\Omega_1, \tilde g_1) \ge \lambda(\Omega_0, g_0) - 2 A
e^{-a l}
\]for  constants $a, A>0$. By (\ref{Ae-a<1/k2}) with $k=0$, this leads to
\[
\lambda(\Omega_1, \tilde g_1) \ge \lambda(\Omega_0, g_0) - 1.
\]

Taking $p>0$ as a variable, the metrics $g^{(p)}_1$ evolves
smoothly with $p$. Lemma \ref{leLstab} shows that
$\lambda(\Omega_1, g^{(p)}_1)$ is a continuous function of $p$.
Observe that
\[
\lambda(\Omega_1, g^{(0)}_1) = \lambda(\Omega_1, \tilde g_1) \ge
\lambda(\Omega_1, g_1) - 1
\]since $g^{(0)}_1  =\tilde g_1$. By the construction of $g^{(p)}_1$, for $x \in H$ such that $x_3 \in [-1, 1]$,
\[
g^{(p)}_1(x)= \left( \frac{1}{2^p} h^2 g_{S^1 \times S^1} \right) \times
g_{R^1}.
\]By Lemma \ref{leLftube}, we know that
\[
\lambda(\Omega_1, g^{(p)}_1) \le \lambda( \frac{1}{2^p} h^2
(S^1 \times S^1) \times [-1, 1]) \to -\infty, \quad p \to \infty.
\]By mean value theorem, there exists a number $p=p_1>0$ so that (\ref{L=L-1}) holds, proving the claim.
This metric $g^{(p_1)}_1$ is the desired metric $g_1$ for
$\Omega_1$, satisfying (\ref{LOm2}).  This completes the
construction of the  component $(\Omega_1, g_1)$, whose
 composition is being summarized here for clarity.
 \be
 \lab{ome1g1}
 \Omega_1 = Z_1  \cup X \cup H \cup Y_1.
 \ee where
 \[
 \al
 g_1 =
 \begin{cases} \text{the round metric} \quad h^2 g_{S^2} \times
 g_{R^1}, \quad
 &\text{on} \quad Z_1 \cup Y_1\\
 g_0, \quad &\text{on} \quad X\\
 \left( \theta^{p_1} h^2 g_{S^1 \times S^1} \right) \times g_{R^1}, \quad
 &\text{on} \quad H.
 \end{cases}
 \eal
 \]\\

 The shape of $\Omega_1$ is depicted here.

 \psfrag{X}{$X$}
 \psfrag{Z}{$Z_1$}
 \psfrag{Y}{$Y_1$}
 \psfrag{H}{$(H, g_1)$}
 \psfrag{O0}{$\Omega_1$}
 \psfrag{zl}{length $l+1$}
 \psfrag{yl}{length $l+1$}
\centerline{\epsfig{figure=tu0.eps, height=5cm}}

Proceeding inductively, suppose we have constructed
\be
 \lab{omekgk}
 \Omega_k = Z_k  \cup X \cup H \cup Y_k.
 \ee where
 \[
 \al
 g_k =
 \begin{cases} \text{the round metric} \quad h^2 g_{S^2} \times
 g_{R^1}, \quad
 &\text{on} \quad Z_k \cup Y_k\\
 g_0, \quad &\text{on} \quad X\\
 \left( \theta^{p_k} h^2 g_{S^1 \times S^1} \right) \times g_{R^1}, \quad
 &\text{on} \quad H,
 \end{cases}
 \eal
 \] and $Z_k$ and $Y_k$ are round necks of length $l+k$.  Now we move to

\medskip

{\it Step 3.} constructing the component $\Omega_{k+1}$ so that

\be \lab{LOmk+1} \lambda(\Omega_{k+1}, g_{k+1}) =
\lambda(\Omega_k, g_k) -\frac{1}{k^2+1}. \ee

This is similar to Step 2, with some modification of parameters. \\

{\it Step 3.1.} Attach the round neck $h^2 S^2 \times [0, 1]$ to
the left end of $\Omega_k$, forming the round neck $h^2 S^2 \times
[0, l+k+1]$ on the left side, which we call $Z_{k+1}$. Then attach
the round neck $h^2 S^2 \times [0, 1]$ to the right end of
$\Omega_k$, forming the round neck $h^2 S^2 \times [0, l+k+1]$ on
the right side, which we call $Y_{k+1}$. The resulting domain is
called $\Omega_{k+1}$ with inherited metric called $\tilde
g_{k+1}$. i.e.

\[
 \Omega_{k+1} = Z_{k+1}  \cup X \cup H \cup Y_{k+1}.
 \] and
 \[
 \al
\tilde g_{k+1} =
 \begin{cases} \text{the round metric} \quad h^2 g_{S^2} \times
 g_{R^1}, \quad
 &\text{on} \quad Z_{k+1} \cup Y_{k+1}\\
 g_0, \quad &\text{on} \quad X\\
 \left( \theta^{p_k} h^2 g_{S^1 \times S^1}\right) \times g_{R^1}, \quad
 &\text{on} \quad H.
 \end{cases}
 \eal
 \]\\

{\it Step 3.2.} Modify $\tilde g_{k+1}$ to a new metric $g_{k+1}$ so that
(\ref{LOmk+1}) holds.

This is again done by pinching $H$, the top portion of the flat
tube. Here are the details.

 Let $\theta$ be the smooth function as in Step 2.  Now consider the metrics on $\Omega_{k+1}$:
\[
\al
g^{(p)}_{k+1}(x)=
\begin{cases}
\tilde g_{k+1}(x), \qquad & x \in \Omega_{k+1} - H\\
\left( \theta^p(x)  h^2 (S^1 \times S^1) \right) \times
g_{R^1},  \qquad & x \in H.
\end{cases}
\eal
\]

We claim that there exists a number $p_{k+1}>0$  so that
 \be
\lab{L=L-1/k2} \lambda(\Omega_{k+1}, g^{(p_{k+1})}_{k+1}) =
\lambda(\Omega_k, g_k) -\frac{1}{k^2+1}.
\ee  Here is the proof.
Regarding $(\Omega_k, g_k)$ as a domain in $(\Omega_{k+1}, \tilde
g_{k+1})$ and applying Lemma \ref{leLX-LY} twice, we know that
\[
\lambda(\Omega_{k+1}, \tilde g_{k+1}) \ge \lambda(\Omega_k, g_k) - 2 A e^{-a (l+k)}
\]for constants $a, A>0$.  Note the length of $Z_k$ and $Y_k$ are $k+l$, which explains
the appearance of the exponential term $e^{-a (l+k)}$.
By (\ref{Ae-a<1/k2}), this leads to
\[
\lambda(\Omega_{k+1}, \tilde g_{k+1}) \ge \lambda(\Omega_k, g_k) -
\frac{1}{k^2+1}.
\]

Taking $p>0$ as a variable, the metrics $g^{(p)}_{k+1}$ evolves smoothly with $p$. Lemma \ref{leLstab}
shows that $\lambda(\Omega_{k+1}, g^{(p)}_{k+1})$ is a continuous function of $p$. Observe that
\[
\lambda(\Omega_{k+1}, g^{(p_k)}_{k+1}) = \lambda(\Omega_{k+1},
\tilde g_{k+1}) \ge \lambda(\Omega_k, g_k) - \frac{1}{k^2+1}
\]since $g^{(p_k)}_{k+1}  =\tilde g_{k+1}$. By the construction of $g^{(p)}_{k+1}$, for $x \in H$ such that $x_3 \in [-1, 1]$,
\[
g^{(p)}_{k+1}(x)=\frac{1}{2^p} h^2 g_{S^1 \times S^1} \times g_{R^1}.
\]From Lemma \ref{leLftube}, we know that
\[
\lambda(\Omega_{k+1}, g^{(p)}_{k+1}) \le \lambda( \frac{1}{2^p}
h^2 (S^1 \times S^1) \times [-1, 1]) \to -\infty, \quad p \to
\infty.
\]By mean value theorem, there exists a number $p=p_{k+1} \ge p_k $ so that (\ref{L=L-1/k2}) holds, proving the claim.
This metric $g^{(p_{k+1})}_{k+1}$ is the desired metric $g_{k+1}$
for $\Omega_{k+1}$, satisfying (\ref{LOmk+1}).  This completes the
construction of the component $(\Omega_{k+1}, g_{k+1})$, finishing
the induction. To summarize,
\be
 \lab{omek+1}
 \Omega_{k+1} = Z_{k+1}  \cup X \cup H \cup Y_{k+1}.
 \ee and
 \[
 \al
 g_{k+1} =
 \begin{cases} \text{the round metric} \quad h^2 g_{S^2} \times
 g_{R^1}, \quad
 &\text{on} \quad Z_{k+1} \cup Y_{k+1}\\
 g_0, \quad &\text{on} \quad X\\
 \left( \theta^{p_{k+1}} h^2 g_{S^1 \times S^1}\right) \times g_{R^1}, \quad
 &\text{on} \quad H.
 \end{cases}
 \eal
 \]

 The shape of $\Omega_{k+1}$ is depicted here.

 \psfrag{X}{$X$}
 \psfrag{Z}{$Z_{k+1}$}
 \psfrag{Y}{$Y_{k+1}$}
 \psfrag{H}{$(H, g_{k+1})$}
 \psfrag{O0}{$\Omega_{k+1}$}
 \psfrag{zl}{length $l+k+1$}
 \psfrag{yl}{length $l+k+1$}
\centerline{\epsfig{figure=tu0.eps, height=5cm}}
\medskip

{\it Step 4.}  pasting together the components to form the manifold $\M$.
See the figure at the end of the step.\\

 In the last step, we have constructed the manifolds $(\Omega_k, g_k)$ for $k=0, 1, 2, 3, ...$.
 Now we define,
 \[
(\Omega_{-k}, g_{-k}) = (\Omega_k, g_k), \quad k=1, 2, ....
 \]
 Finally, we take
\be
\lab{M=Uomgk} \M  = \cup^{\infty}_{k=-\infty} \Omega_k
 \ee
which is the connected, non-overlapping union of  $\Omega_k$, for
all integers $k$ in the
 following pattern.
We connect $\Omega_k$ with $\Omega_{k+1}$ by pasting the right end
of $Y_k$ with left end of $Z_{k+1}$. Here $k=..., -2, -1, 0, 1, 2,
...$. The metric on $\M$, which is inherited from $g_k$, is
denoted by $g$.

It is clear that $\M$ is a complete, connected manifold. Now let us prove $\M$ has
bounded geometry.
Note that except for the top portions of $\Omega_k$, which is denoted  by $H$,  the manifold  $\M$
is consisted of round necks or flat tubes of fixed aperture.
Hence we just need to prove that
$(H, g)$ has bounded geometry.
 The metric $g$ on  $H \subset \Omega_k$ is given by  $(\theta^{p_k}(x) h^2 g_{S^1 \times S^1}) \times g_{R^1}$, where
$\theta(x) = 1/2$ when $x_3 \in [-1, 1]$ and $1/2 \le \theta \le
1$. Write \[
 \lambda_k=\lambda(\Omega_k, g), \quad k=0, 1, 2, ....
 \] Recall by
construction that
\[
\lambda_0 - \Sigma^{|k|}_{j=1} \frac{1}{j^2}  = \lambda_k \le 0,
\quad |k|= 1, 2, ...,
\] which implies
\[
\lambda \left( (1/2)^{p_k} h^2 (S^1 \times S^1) \times [-1, 1] \right) \ge
 \lambda_k \ge \lambda_0 - 10.
\] If $\{ p_k \}$ is unbounded, then by Lemma \ref{leLftube}, the left hand side of the above inequality
tends to $-\infty$ when $k \to \infty$, which leads to a contradiction. Hence $\{ p_k \}$ is a bounded sequence of positive numbers. Since $\theta$ is a smooth bounded function, we know
 $\theta^{p_k}$ has uniformly bounded $C^\infty$ norm. Therefore we have proven that $\M$ has bounded geometry
everywhere. \\
The shape of $\M$ is depicted here.

 \psfrag{X}{$X$}
 \psfrag{Z1}{$Z_k$}
 \psfrag{Y1}{$Y_k$}
 \psfrag{H1}{$(H, g_k)$}
 \psfrag{O1}{$\Omega_k$}
 \psfrag{Z2}{$Z_{k+1}$}
 \psfrag{Y2}{$Y_{k+1}$}
 \psfrag{H2}{$(H, g_{k+1})$}
 \psfrag{O2}{$\Omega_{k+1}$}
 \psfrag{zl}{length $l+1$}
 \psfrag{yl}{length $l+1$}
\centerline{\epsfig{figure=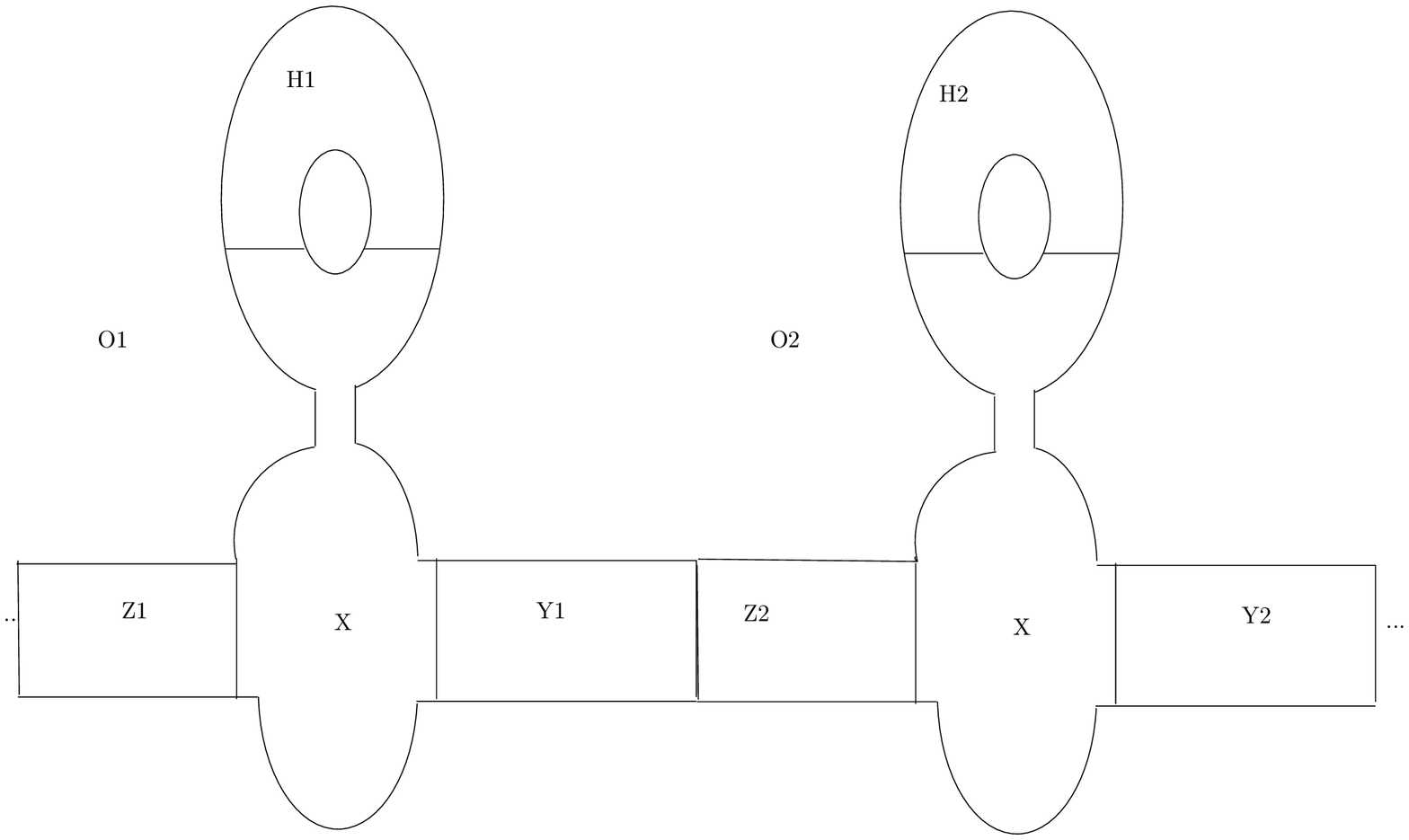, height=5cm, width=14cm}}

{\it Step 5.}  proving that the Log Sobolev functional on $\M$ does not have an extremal. \\

We use the method of contradiction. Suppose that a smooth function
$v$, $\Vert v \Vert_{L^2(\M)}=1$, is an extremal for the Log
Sobolev functional whose infimum is $\lambda=\lambda({\M}, g)$.
Then \[ \lambda = \int_{\M} ( 4 | \nabla v|^2 + R v^2 - v^2 \ln
v^2 ) dg
\] and $v$ is a smooth solution to equation (\ref{eqln}) i.e.
\[
4 \Delta v - R v + 2 v \ln v + \lambda v =0.
\]

Let us recall that
\[
\Omega_k = Z_k \cup X \cup H \cup Y_k,
\]where $Z_k= h^2 S^2 \times [0, l+k]$ and $Y_k= h^2 S^2 \times [0,
l+k]$ are round necks on the left and right side of the core $X$
respectively. In order to distinguish these two necks, we use $z$
to denote points in $Z_k$ with a coordinate $z=(z_1, z_2, z_3)$
described in Definition \ref{decyltub}; and likewise we use $y$ to
denote points in $Y_k$ with a coordinate $y=(y_1, y_2, y_3)$
described in Definition \ref{decyltub}. These two coordinates are
regarded as independent ones.

For each $k=1, 2, ...$, we construct a cut-off function $\eta_k
\in W^{1, \infty}_0(\Omega_k)$ as follows.
\be
\lab{etak=}
\al \eta_k = \begin{cases} \eta_k(z) = z_3, \quad & z \in Z_k, \quad 0
\le z_3 \le 1,\\
\eta_k(z)=1, \quad & z \in Z_k, \quad  1 \le z_3 \le l+k,\\
\eta_k(x)=1, \quad & x \in X \cup H,\\
\eta_k(y) = 1, \quad & y \in Y_k, \quad 0
\le y_3 \le l+k-1,\\
\eta_k(y) = 1-(y_3-l-k+1), \quad & y \in Y_k, \quad  l+k-1 \le y_3
\le l+k.
\end{cases}
\eal \ee The following figure depicts the definition of $\eta_k$.

\psfrag{X}{$X$}
 \psfrag{Z}{$Z_k$}
 \psfrag{Y}{$Y_k$}
 \psfrag{H}{$(H, g_k)$}
 \psfrag{O0}{$\Omega_k$}
 \psfrag{e0}{$\eta_k \le 1$}
 \psfrag{e1}{$\eta_k=1$}
 \psfrag{Zk1}{$Z_{k1}$}
 \psfrag{Yk1}{$Y_{k1}$}
\centerline{\epsfig{figure=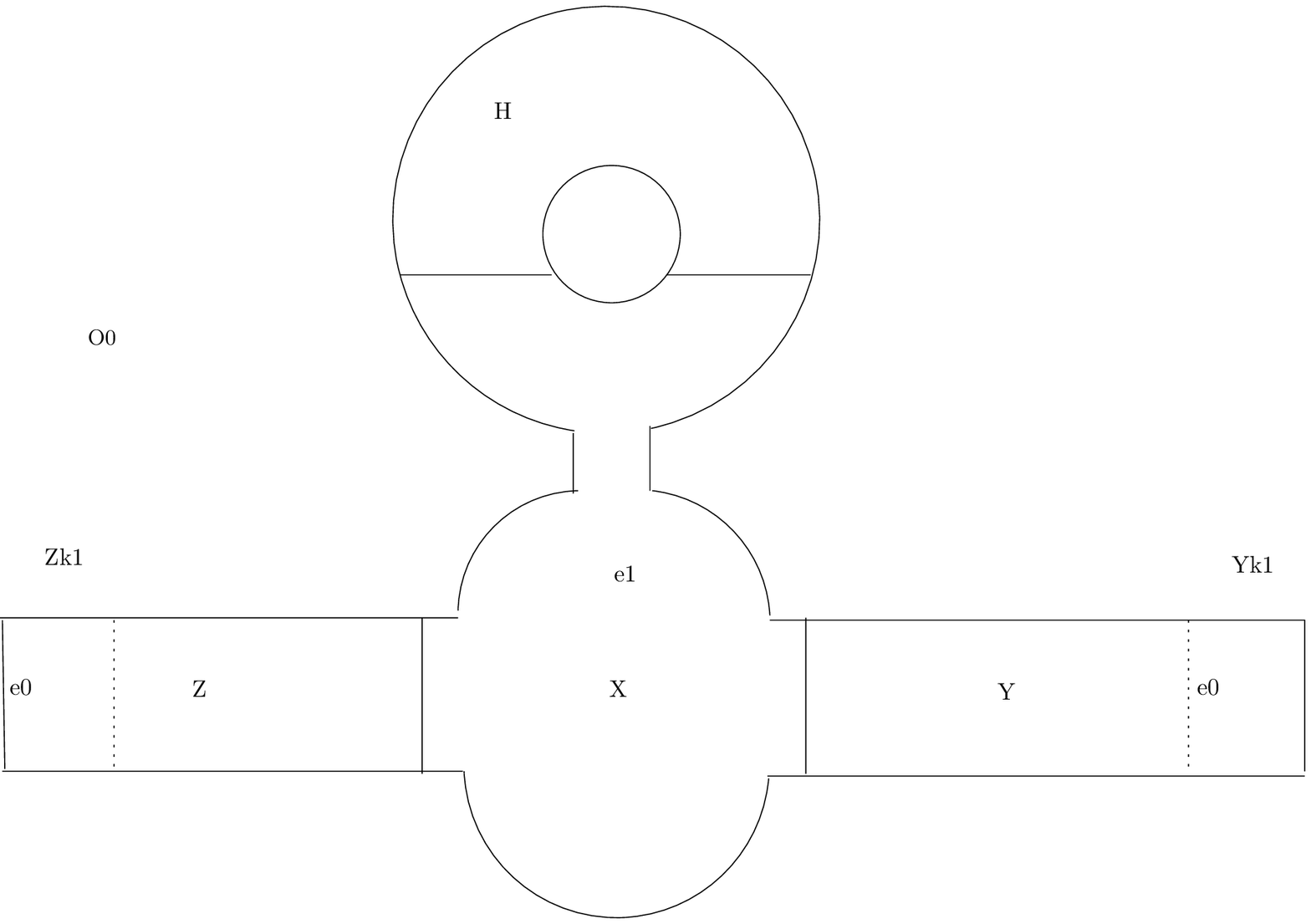, height=5cm, width=14cm}}

\medskip

Since $v$ solves (\ref{eqln}), we can apply Lemma \ref{leLXY} by
taking $E=\Omega_k$ and $F=\M$ there to get
\[
\lambda_k  \int (v \eta_k)^2 dg  \le \lambda  \int (v \eta_k)^2 dg+ 4 \int  v^2 |\nabla \eta_k|^2 dg -
\int (v \eta_k)^2 \ln \eta^2_k dg.
\] Here $\lambda_k = \lambda (\Omega_k, g_k)$. Observe that $|\nabla \eta_k| \le 1$ and that
the function $(\eta_k)^2 \ln \eta^2_k$, which is nonzero only in
the support of $\nabla \eta_k$, is bounded from below by  $-
e^{-1}$. Therefore
\be
\lab{Lk<L} \lambda_k  \int (v \eta_k)^2 dg
\le \lambda  \int (v \eta_k)^2 dg+ 5 \int_{supp \nabla \eta_k} v^2
dg.
\ee By definition of $\eta_k$,  $supp \nabla \eta_k$ is the
disjoint union of two short round necks, i.e.
\be
\lab{suppetak} supp \nabla \eta_k = Z_{k1} \cup Y_{k1} \ee  where
\[
Z_{k1} \equiv  \{ z \in Z_k \, |  0 \le z_3 \le 1\}, \qquad Y_{k1}
\equiv  \{ y \in Y_k \, |  l+k-1 \le y_3 \le l+k \}.
\]  Hence (\ref{Lk<L}) implies
\be \lab{Lk-L<} (\lambda_k -\lambda)  \int (v \eta_k)^2 dg  \le 5
\int_{Z_{k1} }  v^2 dg + 5 \int_{Y_{k1} }  v^2 dg = 5 \int_{supp
\nabla \eta_k } v^2 dg. \ee

Next we prove that the right hand side of (\ref{Lk-L<}) is
exponentially small. Observe that $Z_{k1}$ is a middle segment of
$Y_{k-1} \cup Z_k$, which is, when writing in one coordinate, a
round neck of the form $h^2 S^2 \times [0, 2l+2k-1]$.  The
segments
\[
W_{k-1} \equiv \{y \in Y_{k-1} | 0 \le y_3 \le 2\} \quad
\text{and} \quad
 E_k \equiv \{z \in Z_k | k+l-2 \le z_3 \le k+l\}
\]are at the left and right end of the round neck
  respectively.  See the figure below. \\

 \psfrag{Y}{$Y_{k-1}$}
 \psfrag{Z}{$Z_k$}
 \psfrag{YZ}{$Y_{k-1} \cup Z_k$}
 \psfrag{Zk1}{$Z_{k1}$}
 \psfrag{y3}{$W_{k-1}$}
 \psfrag{z3}{$E_k$}
 \psfrag{O1}{$\Omega_{k-1}$}
 \psfrag{O2}{$\Omega_k$}
 \centerline{\epsfig{figure=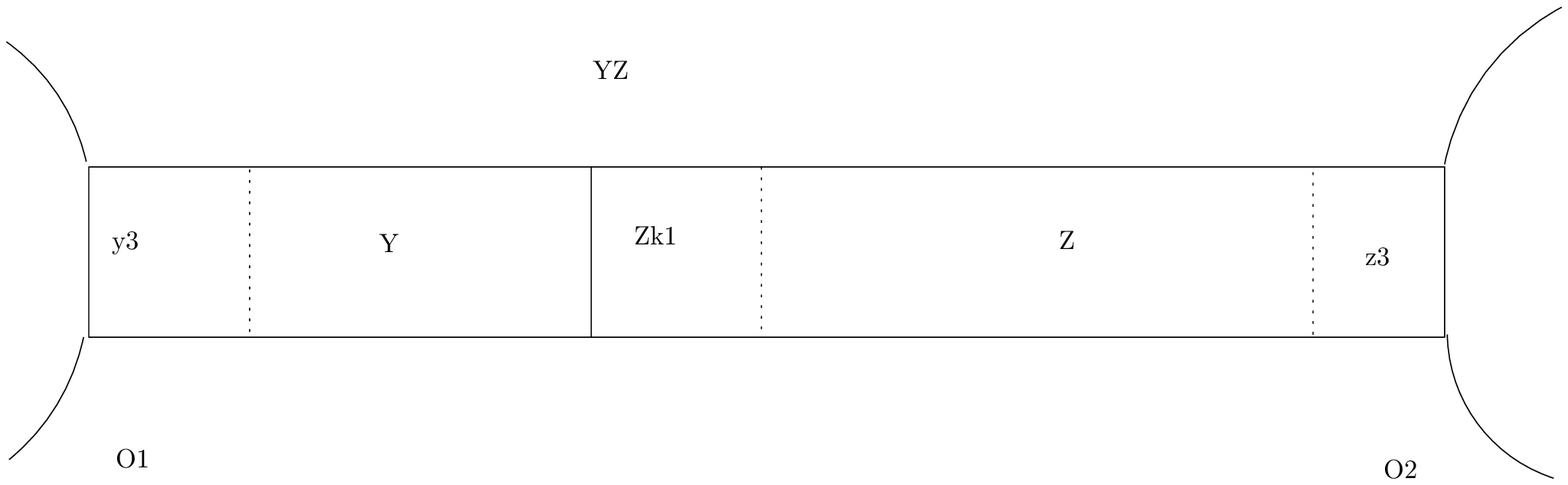, height=4cm, width=12cm}}
 \medskip

By Lemma \ref{lev<e-l}, we have
\[
\int_{Z_{k1}} v^2 dg \le A e^{-a (l+k-1)} \left[ \int_{ \{y \in Y_{k-1} | 0 \le y_3 \le 2\}} v^2 dg +
\int_{ \{z \in Z_k | k+l-2 \le z_3 \le k+l\}} v^2 dg \right].
\] Note from (\ref{etak=}) that
\[
\eta_{k-1}=1 \quad \text{in} \quad  W_{k-1}=\{y \in Y_{k-1} | 0
\le y_3 \le 2\} \subset Y_{k-1} \subset \Omega_{k-1},
\]
\[
\eta_k=1 \quad  \text{in} \quad
 E_k =\{z \in Z_k | k+l-1 \le z_3 \le k+l\} \subset Z_k \subset \Omega_k.
\]Hence
\be
\lab{intZv2}
\int_{Z_{k1}} v^2 dg \le A e^{-a (l+k-1)} \left[ \int_{ \Omega_{k-1}} (\eta_{k-1} v)^2 dg +
\int_{ \Omega_k} (\eta_k v)^2 dg \right].
\ee

Similarly, we see that $Y_{k1}$ is a middle segment of $Y_k \cup Z_{k+1}$, which is, when writing in one coordinate, a round neck of the form
$h^2 S^2 \times [0, 2l+2k+1]$.  The segments
\[
 W_k \equiv \{y \in Y_k | 0 \le y_3 \le 2\} \quad \text{and} \quad
 E_{k+1} \equiv \{z \in Z_{k+1} | k+l-1 \le z_3 \le k+l+1\}
\] are  the left and right end of the round neck.
See the figure below. \\

 \psfrag{Y}{$Y_k$}
 \psfrag{Z}{$Z_{k+1}$}
 \psfrag{YZ}{$Y_k \cup Z_{k+1}$}
 \psfrag{Yk1}{$Y_{k1}$}
 \psfrag{y3}{$W_k$}
 \psfrag{z3}{$E_{k+1}$}
 \psfrag{O1}{$\Omega_k$}
 \psfrag{O2}{$\Omega_{k+1}$}
 \centerline{\epsfig{figure=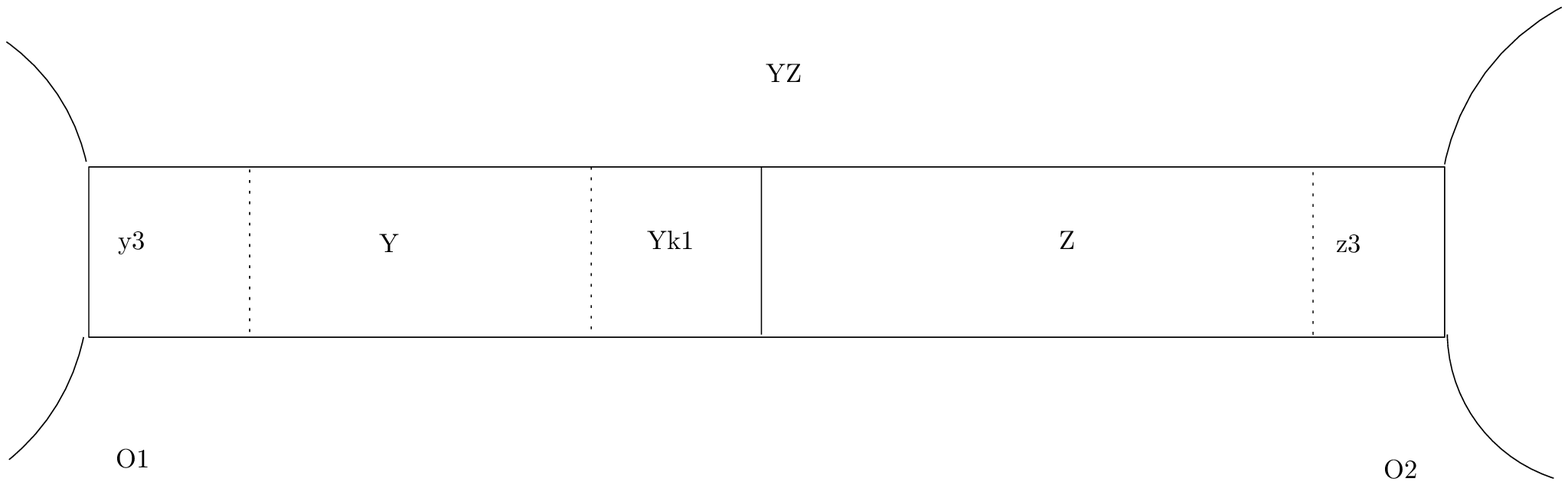, height=4cm, width=12cm}}
 \medskip

 By Lemma \ref{lev<e-l}, we have
\[
\int_{Y_{k1}} v^2 dg \le A e^{-a (l+k)} \left[ \int_{ \{y \in Y_k | 0 \le y_3 \le 2\}} v^2 dg +
\int_{ \{z \in Z_{k+1} | k+l-1 \le z_3 \le k+l+1\}} v^2 dg \right].
\] Note that
\[
\eta_k=1 \quad \text{in} \quad  W_k= \{y \in Y_k | 0 \le y_3 \le
2\} \subset Y_k \subset \Omega_k
\] and
\[
\eta_{k+1}=1  \quad \text{in} \quad
 E_{k+1}=\{z \in Z_{k+1} | k+l-1 \le z_3 \le k+l+1\} \subset Z_{k+1} \subset \Omega_{k+1}.
 \] Hence
  \be
\lab{intYv2} \int_{Y_{k1}} v^2 dg \le A e^{-a (l+k)} \left[ \int_{
\Omega_k} (\eta_k v)^2 dg + \int_{ \Omega_{k+1}} (\eta_{k+1} v)^2
dg \right]. \ee  By this, (\ref{intZv2}) and (\ref{suppetak}), we
obtain
\[
\int_{supp \nabla \eta_k} v^2 dg \le A e^{-a (l+k-1)} \left[
\int_{ \Omega_{k-1}} (\eta_{k-1} v)^2 dg + 2 \int_{\Omega_k}
(\eta_k v)^2 dg + \int_{ \Omega_{k+1}} (\eta_{k+1} v)^2 dg\right]
\] where $k=1, 2, 3, ...$. Recall that $(\Omega_{-k}, g_{-k}) = (\Omega_k,
g_k)$ be definition. Therefore, we can derive, in a similar
manner,
\[
\int_{supp \nabla \eta_k} v^2 dg \le A e^{-a (l+|k|-1)} \left[
\int_{ \Omega_{k-1}} (\eta_{k-1} v)^2 dg + 2 \int_{\Omega_k}
(\eta_k v)^2 dg + \int_{ \Omega_{k+1}} (\eta_{k+1} v)^2 dg\right]
\] where $k=0, -1, -2, -3, ...$. Adding the last two inequalities
together, we deduce
\[
\al
 \Sigma^\infty_{k=-\infty} \int_{supp \nabla \eta_k}  v^2 dg &=
  \Sigma^\infty_{k=-\infty} \int_{Z_{k1}}  v^2 dg + \Sigma^\infty_{k=-\infty} \int_{Y_{k1}}  v^2 dg \\
&   \le
 4 A  e^{2 a} \Sigma^\infty_{k=-\infty} e^{- a (l+|k|)}
\int_{\Omega_k} (v \eta_k)^2 dg.
\eal
\]
By (\ref{Lk-L<}), this implies
\[
 \Sigma^\infty_{k=-\infty} (\lambda_k -\lambda - 20 A e^{2a} e^{- a (l+k)}  )  \int (v \eta_k)^2 dg \le 0.
\]Recall, by construction,
\[
\lambda_k -\lambda \ge \lambda_k -\lambda_{k+1} = \frac{1}{k^2+1},
\quad k=1, 2, 3, ...
\]and
\[
\lambda_k -\lambda \ge \lambda_k -\lambda_{k-1} = \frac{1}{k^2+1},
\quad k=-1, -2, -3, ... ,
\] and
\[
\lambda_0 -\lambda  \ge \lambda_0 - \lambda_1 =1.
\]So finally
we deduce
\[
 \Sigma^\infty_{k=-\infty} ( \frac{1}{1+k^2} - 20 A e^{2 a} e^{- a (l+|k|)}  )  \int (v \eta_k)^2 dg \le 0.
\] This is a contradiction because $\frac{1}{1+k^2} - 20 A e^{2a} e^{- a (l+|k|)}>0$ by our choice of $l$ in
(\ref{Ae-a<1/k2}). Therefore no extremal for the Log Sobolev
functional exists. \qed

\section{$W$ entropy and a no breather result for noncompact Ricci flow}

In this section we discuss some applications of Theorem \ref{thmain} to
Perelman's $W$ entropy and Hamilton's Ricci flow. We will use the following notations.
 $g=g(t)$ is a metric which evolves with time; $d(x, y, t)$ or $d(x, y, g(t))$ will denote the
corresponding distance function; $dg(t)$ denotes the volume
element under $g(t)$; We will still use $\nabla$, $\Delta$ the
corresponding gradient and Laplace-Beltrami operator, when no
confusion arises.

The following definition is one of several equivalent ways in
which Perelman's W entropy can be written.

\begin{definition}  (W  entropy)
\lab{deW}
Let $v \in W^{1, 2}(\M)$ and $\tau>0$ be a parameter. The $W$ entropy is the quantity
\be
\lab{Wshang}
W(g, v, \tau) \equiv  \int_{\bf M}  \left[ \tau (4 |\nabla v|^2 + R v^2) - v^2 \ln v^2
-\frac{n}{2} (\ln 4 \pi \tau) \ v^2  -n v^2 \right]
dg.
\ee
\end{definition}

Let $c>0$ be a positive constant, it is clear that the W entropy has the following scaling invariant
property
\[
W(c g, c^{-n/2} v, c \tau) = W(g, v, \tau).
\] Hence we can always take $\tau=1$ if necessary.  If $\tau=1$ and $\Vert v \Vert_{L^2(\M)} = 1$, then
\be
\al
\lab{W=L}
W(g, v, 1) &=  \int_{\bf M}  \left[  (4 |\nabla v|^2 + R v^2) - v^2 \ln v^2
\right] dg -\frac{n}{2} (\ln 4 \pi)  -n \\
&= L(v, g) - \frac{n}{2} (\ln 4 \pi)  -n.
\eal
\ee Here $L(v, g)$ is the Log Sobolev functional given in (\ref{logfanhan}).
Therefore, the W entropy and the Log Sobolev functional differ only by a normalizing constant after
scaling.

Perelman also introduced the so called $\mu$ invariant.

\begin{definition}
Given a noncompact manifold $(\M, g)$ and parameter $\tau>0$,
the $\mu$ invariant is the quantity
\[
\mu(g, \tau)\\
 = \inf \{  W(g, v, \tau)  \, | \,
v \in C^\infty_0(\M), \, \Vert v \Vert_{L^2(\M)=1} \}.
\]
\end{definition}

In view of Definition \ref{delamb}, we introduce
 $\mu$ invariant near infinity.

\begin{definition}

Given a noncompact manifold $(\M, g)$ and parameter $\tau>0$, the $\mu$ invariant at infinity is the quantity
\[
\al
\mu_\infty(g, \tau)
= &\lim_{r \to \infty}
 \inf \{  \int_{\M-B(0, r)}
  \left[ \tau (4 |\nabla v|^2 + R v^2) - v^2 \ln v^2
-\frac{n}{2} (\ln 4 \pi \tau) \ v^2  -n v^2 \right]
dg  \,  | \\
 & \qquad \quad  \quad
v \in C^\infty_0(\M-B(0, r)), \quad \Vert v \Vert_{L^2(\M-B(0, r))=1} \}.
\eal
\]
\end{definition}

Since the $W$ entropy and the Log Sobolev functional differ only by a constant after scaling, Theorem
\ref{thmain} can be immediately transplanted as

\begin{theorem}
\lab{thWshang}
  (a). Let $\M$ be a complete, connected
noncompact manifold with bounded geometry, and $\tau>0$ be a parameter.
 Suppose $\mu({\M}, \tau)<\lambda_\infty({\M}, \tau)$,
then there exists a smooth extremal $v$ for the $W$ entropy (\ref{Wshang}).
Also, there exist positive constants $a, A>0$ and a point $0 \in \M$ such that
\[
v(x) \le A e^{- a d^2(x, 0)}.
\]

(b).  There exists  a complete, connected noncompact manifold with
bounded geometry such that $\mu({\M}, \tau)<\lambda_\infty({\M},
\tau)$, but the $W$ entropy (\ref{Wshang}) does not have an
extremal.
\end{theorem}

In the rest of the section, we describe two more applications of this theorem.
The first one is an extension of Perelman's monotonicity formula for the $W$ entropy from
the compact case to some noncompact ones.

Let us briefly recall Perelman's monotonicity formula.  Consider the final value problem of the conjugate heat equation
 coupled with the Ricci flow $({\M}, g(t))$ on a compact manifold $\M$ and on the time interval
 $[t_1, t_2]$.
\be
\lab{conjheat}
\begin{cases}
\Delta u - R u + u_t=0, \quad t \in [t_1,  t_2]\\
u(x, t_2) = u_2\\
\partial_t g(t) = - 2 Ric, \quad t \in [t_1, t_2].
\end{cases}
\ee Here $\Delta$ is the Laplace-Beltrami operator with respect to the metric $g(t)$;  $R$
and $Ric$ are the scalar curvature  and Ricci curvature with respect to $g(t)$; and $u_2=u_2(x)$ is a smooth
function such that $\Vert u_2 \Vert_{L^1({\M}, g(t_2))}=1$.
In the definition of the $W$ entropy, we take $\tau=L-t$ and $v(\cdot, t)=\sqrt{u(\cdot, t)}$.
Perelman (\cite{P:1} section 3) proved that
\be
\lab{dtW=}
\frac{d}{dt} W(g(t), v(\cdot, t), L-t)=
2 \tau \int_{\bf M} \left| Ric -Hess \ln u - \frac{1}{2 \tau} g \right|^2
\ u \  dg(t).
\ee

If $\M$ is noncompact, then the above formula needs certain
justification. One reason is that the term $Hess \ln u$ may grow
to infinity and hence the integral may diverge. Consequently,
certain extra decay conditions are needed on $u$ and $Hess \ln u$.
When $u$ is the fundamental solution of the conjugate heat
equation, a noncompact version of the above formula has been
carefully established in \cite{CCGGIIKLLN3:1} Chapters 19, 20, 21
and the paper \cite{CTY:1}. They employed  a number of technical
tools such as Log gradient bounds for positive solutions of
(\ref{conjheat}) and pointwise bounds on the fundamental solution
of (\ref{conjheat}). With the help of these tools and the decay
estimate of extremals of the $W$ entropy, we extend (\ref{dtW=})
to a noncompact case where the final value $u_2$ is the square of
an extremal of the $W$ entropy. The point of the following
corollary is that once an extremal exists, then no other decay
conditions are needed.

\begin{corollary}
\lab{coclaim} Let $({\M}, g(t))$ be a Ricci flow which has bounded
geometry in the finite time interval $[t_1, t_2]$. Assume also
that the 4-th order derivatives of the curvature tensor are
uniformly bounded in ${\M} \times [t_1, t_2]$.
 Let $\tau=L-t$ with $L>t_2$
be a parameter. Suppose the $W$ entropy $W(g(t_2), v, T-t_2)$ has
an extremal $v_2$.  Let $u$ be the solution of the final value
problem of the conjugate heat equation:
\[
\begin{cases}
\Delta u - R u + u_t=0, \quad t \in [t_1, t_2]\\
u(x, t_2) = v^2_2\\
\partial_t g(t) = - 2 Ric, \quad t \in [t_1, t_2].
\end{cases}
\] Let $v=v(x, t) = \sqrt{u(x, t)}$. Then, for all $t \in [t_1,
t_2]$, the $W$ entropy $W(g(t), v, T-t)$ is well defined. Moreover
\[
\frac{d}{dt} W(g(t), v, L-t) =2 \tau \int_{\bf M} \left| Ric -Hess
\ln u - \frac{1}{2 \tau} g \right|^2 \ u \  dg(t). \]
 \proof
\end{corollary}

 The task is to show that relevant integrands has quadratic exponential decay at infinity.
 After this, the proof is the same as Perelman's in the compact
 case.

{\it Step 1.}  First we show that there exist positive constants $A_1, a_1$ and a point $0 \in \M$ such that
\be
\lab{ubounds}
  u(x, t) \le A_1 e^{-a_1 d^2(x, 0, t)}.
 \ee  This bound follows from the decay of the extremal $v_2$ in Theorem \ref{thWshang} (a) and
 the following bounds on $G=G(x, t; y, t_2)$,  the fundamental solution of the conjugate heat equation (\ref{conjheat}).
 Observe that the Ricci flow has bounded geometry in the finite time interval $[t_1, t_2]$. Hence
 the distance functions $d(x, 0, t)$ are equivalent when $t \in [t_1, t_2]$. The same can be said for
 volumes $|B(x, r, t)|_{g(t)}$. By
 \cite{CCGGIIKLLN3:1} Chapters 19
or  \cite{CTY:1} Section 5, there are the bounds:
\[
\al
G(x, t; y, t_2) &\ge  \frac{1}{\alpha \sqrt{|B(x, \sqrt{t_2-t}, t)|_{g(t)}}  \sqrt{|B(y, \sqrt{t_2-t}, t)|_{g(t)}}} e^{- \frac{d^2(x, y, t)}{\beta (t_2-t)}},\\
 G(x, t; y, t_2) &\le \frac{\alpha}{\sqrt{|B(x, \sqrt{t_2-t}, t)|_{g(t)}}  \sqrt{|B(y, \sqrt{t_2-t}, t)|_{g(t)}}} e^{-\beta \frac{d^2(x, y, t)}{(t_2-t)}},
\eal
\] where the constants $\alpha$ and $\beta$ depend on $\M$, $t_1$ and $t_2$. These bounds can be
regarded as generalization of the bounds in \cite{LY:1} for the heat equation under fixed metrics.
By the assumption of bounded geometry and classical volume comparison theorem, there
exist positive constants $c$,  $c_1$ and $c_2$ such that
\[
c_1 \min \{1, (t_2-t)^{n/2} \} \le |B(x, \sqrt{t_2-t}, t)|_{g(t)} \le c (t_2-t)^{n/2} e^{c_2  \sqrt{t_2-t}},
\]
\[
c_1 \min \{1, (t_2-t)^{n/2} \} \le |B(y, \sqrt{t_2-t}, t)|_{g(t)} \le c (t_2-t)^{n/2} e^{c_2  \sqrt{t_2-t}},
\] Hence we have the bounds: for $t \in [t_1, t_2]$ and $x, y \in \M$,
\be
\lab{Gbounds}
 \frac{1}{\alpha (t_2-t)^{n/2}} e^{-\beta \frac{d^2(x, y, t)}{(t_2-t)}} \le G(x, t; y, t_2) \le \frac{\alpha}{(t_2-t)^{n/2}} e^{-\beta \frac{d^2(x, y, t)}{(t_2-t)}},
\ee  where the constant $\alpha=\alpha({\M}, t_1, t_2)$ may have changed from its previous value.
Therefore,
\[
u(x, t) =\int_{\M} G(x, t; y, t_2) u(y, t_2) dg(t_2) \le
\int_{\M} \frac{\alpha}{(t_2-t)^{n/2}} e^{-\beta \frac{d^2(x, y, t)}{(t_2-t)}} u(y, t_2) dg(t_2)
\]

By Theorem \ref{thWshang} (a) (in fact Lemma \ref{ledecay} is sufficient),
 there exist positive constants $a, A>0$  such that
\be
\lab{v2decay}
u(x, t_2) = v^2_2(x) \le  2 A e^{- 2 a d^2(x, 0, t_2)}.
\ee  The last two inequalities imply
\be \lab{uxt<} u(x, t)  \le 2 \alpha A \int_{\M}
\frac{1}{(t_2-t)^{n/2}} e^{-\beta \frac{d^2(x, y, t)}{(t_2-t)}}
e^{- 2 a d^2(y, 0, t_2)} dg(t_2). \ee By triangle inequality,
there exist $a_1>0$ such that
\[
-\beta \frac{d^2(x, y, t)}{(t_2-t)} - 2 a d^2(y, 0, t_2) \le -a_1 d^2(x, 0, t_2)
-\beta \frac{d^2(x, y, t)}{2(t_2-t)}.
\] Here we used the fact that distances at different time levels are comparable again.
Hence
\[
\al &\int_{\M} \frac{1}{(t_2-t)^{n/2}} e^{-\beta \frac{d^2(x, y,
t)}{(t_2-t)}} e^{- 2 a d^2(y, 0, t_2)} dg(t_2) \\
&\le e^{ -a_1 d^2(x, 0, t_2) } \int_{\M} \frac{1}{(t_2-t)^{n/2}}
e^{-\beta \frac{d^2(x, y, t)}{2(t_2-t)}} dg(t_2)\\
&=e^{ -a_1 d^2(x, 0, t_2) } [ \Sigma^\infty_{k=0} \int_{ 2^{k-1}
\sqrt{t_2-t} \le d(x, y, t) \le 2^k \sqrt{t_2-t}}
\frac{1}{(t_2-t)^{n/2}} e^{-\beta \frac{d^2(x, y, t)}{2(t_2-t)}}
dg(t_2) \\
&\qquad \qquad \qquad + \int_{ d(x, y, t) \le  \sqrt{t_2-t}}
\frac{1}{(t_2-t)^{n/2}} e^{-\beta \frac{d^2(x, y, t)}{2(t_2-t)}}
dg(t_2) ].
 \eal
\] Since $\M$ has bounded geometry, the classical volume
comparison theorem tells us
\[
|B(x, 2^k \sqrt{t_2-t}, t)|_{g(t_2)} \le C e^{c 2^k \sqrt{t_2-t}}
(t_2-t)^{n/2}.
\] Here we just used the fact that volume elements at different
time levels in $[t_1, t_2]$ are equivalent. Hence
\[
\al &\int_{\M} \frac{1}{(t_2-t)^{n/2}} e^{-\beta \frac{d^2(x, y,
t)}{(t_2-t)}} e^{- 2 a d^2(y, 0, t_2)} dg(t_2) \\
&\le e^{ -a_1 d^2(x, 0, t_2) } \left[ \Sigma^\infty_{k=0} C e^{c
2^k \sqrt{t_2-t}}
 e^{-\beta 2^{2(k-1)}/2} + C e^{c \sqrt{t_2-t_1}}\right],
 \eal
\] which shows
\be
\lab{convoldecay}
 \int_{\M} \frac{1}{(t_2-t)^{n/2}} e^{-\beta
\frac{d^2(x, y, t)}{(t_2-t)}} e^{- 2 a d^2(y, 0, t_2)} dg(t_2) \le
C e^{ -a_1 d^2(x, 0, t_2) },
 \ee where $C$ depends on $t_2-t_1$.
 Substituting this to (\ref{uxt<}), we
deduce
\[
u(x, t)  \le A_1 e^{ -a_1 d^2(x, 0, t_2) } .
\]This proves the  bound in (\ref{ubounds}). \\

{\it Step 2.} We prove that the integrand in the $W$ entropy has
quadratic exponential decay.

 For convenience, we denote the integrand in the $W$ entropy as
\be
\lab{iu=}
i(u) =i(u)(x, t) \equiv \left[ \tau ( \frac{|\nabla u|^2}{u} + Ru) - u \ln u - \frac{n}{2} \ln(4 \pi \tau) u - n u \right]
(x, t).
\ee Here we have used the relation that $u=v^2$ on (\ref{deW}).
We now prove that there exist positive constants
$A_1$ and $a_1$ such that
\be
\lab{iu<}
| i(u)(x, t) | \le A_1 e^{-a_1 d^2(x, 0, t)}.
\ee  By the bound (\ref{ubounds}), we know that the term $u \ln u$ satisfies
\[
| u \ln u| =  \sqrt{u} \sqrt{u} | \ln u| \le C \sqrt{u} \le A_1
e^{-a_1 d^2(x, 0, t)},
\]whence it also has quadratic exponential decay. Here the values of $A_1$ and $a_1$ may have changed.
So it suffices to prove that the term $\frac{|\nabla u|^2}{u}$ has the decay too.

To this end, we
recall by direct computation (see  Proposition 6.1.2 in \cite{Z:1} e.g.) that
\[
\al
&H^*(\frac{|\nabla u|^2}{u} + R u)\\
& = \frac{2}{u} \left( u_{ij}-
\frac{u_i u_j}{u} \right)^2 + 2 \nabla R \nabla u + \frac{4}{u}
Ric(\nabla u, \nabla u) + 2 |Ric|^2 u + 2 \nabla R \nabla u + 2 u
\Delta R.
\eal
\] Here $H^* = \Delta - R + \partial_t$ is the conjugate heat operator.
Thus
\[
H^*(\frac{|\nabla u|^2}{u} + R u) \ge  - K_1 (| \nabla u| + |u| +
\frac{|\nabla u|^2}{u}),
\]where the constant $K_1 (\ge 0)$ depends on the supremum of
 $| \nabla R|$, \ $|\Delta R |$ and the lower bound of $Ric$.
Since $| \nabla u| \le \frac{|\nabla u|^2}{u} + u$, we deduce
\be
\lab{H*du2/u}
H^*(\frac{|\nabla u|^2}{u} + R u) \ge  - K_1
(\frac{|\nabla u|^2}{u} + R u) - K_2 u,
\ee  where $K_2$ depends on
$K_1$, the supremum of $|R|$ and $u$. We mention that all the
curvatures involved here are bounded according to our assumption on the Ricci flow.

At time $t_2$,  $u(x, t_2) = v^2_2$.  Hence $\frac{|\nabla u|^2}{u} + R u =
4 | \nabla v_2|^2 + R v^2_2.$  Notice that $v_2$ satisfies
 the equation for extremals: for $\tau = L-t_2$,
\[
\tau ( 4 \Delta v_2 - R v_2) + 2 v_2 \ln v_2 + \frac{n}{2} \ln(4 \pi \tau) v_2 + n v_2 + \mu v_2=0
\] By Lemma \ref{lemvi} part (b), we have
\[
\sup_{B(x, 1/2, t_2)}  | \nabla v_2|^2 \le  C \int_{B(x, 1, t_2)} v^2_2 dg(t_2) \le C A^2 e^{- 2 a d^2(x, 0, t_2)},
\] where the last inequality is due to the decay of $v_2$ in (\ref{v2decay}).
By this and the decay of $v_2$ again, we know that, at time $t_2$,
\be
\lab{Qt2<}
\left| \frac{|\nabla u|^2}{u} + R u \right|(x, t_2)
  \le A_1 e^{- a_1 d^2(x, 0, t_2)}.
\ee Now define
\[
Q=Q(u)= e^{K_1 t} \left( \frac{|\nabla u|^2}{u} + R u \right).
\] By (\ref{H*du2/u}) and (\ref{Qt2<}), we know that
\[
\begin{cases}
\Delta Q - R Q + \partial_t Q \ge  - K_2 e^{K_1t} u, \quad t \in [t_1, t_2],\\
Q(\cdot, t_2) \le  A_1e^{K_1 t_2} e^{- a_1 d^2(x, 0, t_2)}.
\end{cases}
\] By the maximum principle (see \cite{CCGGIIKLLN2:1} Chapter 12 e.g.),
this implies, for $t \in  [t_1, t_2]$,
\be
 \lab{Qxt<}\al &Q(x, t)
\le \int_{\M} G(x, t; y, t_2)  A_1e^{K_1 t_2} e^{- a_1 d^2(y, 0,
t_2)} dg(t_2) \\
&\qquad \qquad + \int^{t_2}_t \int_{\M} G(x, t; y, s)  K_2 e^{K_1
t_2} u(y, s) dg(s) ds. \eal \ee We mention that even though $Q$ is
a smooth function, it may not be a bounded one for each time
level, due to the appearance of the term $\frac{|\nabla u|^2}{u}$.
In order to apply the maximum principle, one needs some growth
condition on $Q$ near infinity. The way to justify (\ref{Qxt<}) is
to replace $u$ by the function $u_\e$ which is the solution to
\[
\begin{cases}
\Delta u_\e - R u_\e +  \partial_t u_\e =0, \quad t \in [t_1, t_2]\\
u_\e(x, t_2) = v^2_2 +\e \\
\partial_t g(t) = - 2 Ric, \quad t \in [t_1, t_2].
\end{cases}
\]  Here $\e>0$ is a positive number. It is clear that $u_\e \to u$ pointwise when $\e \to 0$.
Since $u_\e$ is bounded from above and below by positive constants, we know that
$Q_\e = Q(u_\e)$ is a bounded function. Moreover it holds
\[
\begin{cases}
\Delta Q_\e - R Q_\e + \partial_t Q_\e \ge  - K_2 e^{K_1t} u_\e, \quad t \in [t_1, t_2],\\
Q_\e(\cdot, t_2) \le  A_1e^{K_1 t_2} e^{- a_1 d^2(x, 0, t_2)} + C \e.
\end{cases}
\] Now we can apply the maximum principle for $Q_\e$ to derive
\be
 \lab{Qext<}\al &Q_\e(x, t) \le \int_{\M} G(x, t; y, t_2)
Q_\e(y, t_2) dg(t_2) \\
&\qquad \qquad + \int^{t_2}_t \int_{\M} G(x, t; y, s)  K_2 e^{K_1
t_2} u_\e(y, s) dg(s) ds. \eal \ee Taking $\e \to \infty$, this
implies (\ref{Qxt<}).

By (\ref{Qxt<}),  (\ref{Gbounds}) and (\ref{ubounds}), we derive
\[
\al
Q(x, t) &\le A_1e^{K_1 t_2} \int_{\M} \frac{\alpha}{(t_2-t)^{n/2}} e^{-\beta \frac{d^2(x, y, t)}{(t_2-t)}}   e^{- a_1 d^2(y, 0, t_2)} dg(t_2) \\
&\qquad
+  K_2 e^{K_1 t_2} \int^{t_2}_t \int_{\M} \frac{\alpha}{(s-t)^{n/2}} e^{-\beta \frac{d^2(x, y, t)}{(s-t)}}  A_1 e^{- a_1 d^2(y, 0, s)} dg(s) ds.
\eal
\] Using the fact that distance functions and volume elements at different time levels are equivalent,
we can apply (\ref{convoldecay}) to the above inequality to deduce
\[
Q(x, t) \le C e^{- c d^2(x, 0, t)}
\] where $c$ and $C$ are positive constants which may depend on $t_1$ and $t_2$.
This proves  that
\[
\left| \frac{|\nabla u|^2}{u} + R u \right|(x, t) \le C e^{- c d^2(x, 0, t)}, \quad t \in [t_1, t_2],
\]which implies (\ref{iu<}). \\

{\it Step 3.}  Completion of the proof.

Let $u$ and $\tau$ be the same as in the statement of the Corollary.
In the paper \cite{P:1} Proposition 9.1, Perelman introduced the quantity
\be
\lab{pu=}
P(u) = \tau (- 2 \Delta u + \frac{|\nabla u|^2}{u} + Ru) - u \ln u - \frac{n}{2} \ln(4 \pi \tau) u - n u
\ee and proved that
\be
\lab{H*Pu}
H^* P(u) = 2 \tau \left| Ric - Hess \ln u -\frac{g}{2 \tau} \right|^2 u.
\ee We mention that in \cite{P:1}, the quantity $P(u)$ here is denoted by $v=v(f)$ where $f$ is determined by
$u = \frac{e^{-f}}{(4 \pi \tau)^{n/2}}$.
Observe that
\be
\lab{puiu}
P(u) = - 2 \tau \Delta u + i(u)
\ee where $i(u)$ is the integrand of the $W$ entropy used in the previous step.

Next we will integrate (\ref{H*Pu}). However, at the moment, we do not know the if terms
involved are integrable. So we need to use certain cut off function.
Let $L=L(x)$ be a smooth function on $\M$ such that
\[
\al
&|\nabla L(x) | + |\nabla^2 L(x) | + |\nabla^3 L(x) | + |\nabla^4 L(x) | \le C_1, \qquad x \in \M, \\
&C^{-1}_1 L(x) \le d(x, 0, g(t_1)) \le C_1 L(x), \qquad x \in \M.
\eal
\] Here the covariant derivatives are with respect to $g(t_1)$. Under our assumption of bounded geometry, it is well known that such a function exists. See for example Proposition 19.37 and the remark
right after it in \cite{CCGGIIKLLN3:1}.  By our assumption of uniformly bounded curvature and
its up to  $4$-th order derivatives, it is easy to check that, there exists $C_2>0$ depending on $t_1$ and $t_2$
such that
\be
\lab{LLLL}
\al
&|\nabla L(x) | + |\nabla^2 L(x) | + |\nabla^3 L(x) | + |\nabla^4 L(x) | \le C_2, \qquad x \in \M, \\
&C^{-1}_2 L(x) \le d(x, 0, g(t)) \le C_2 L(x), \qquad x \in \M.
\eal
\ee Here the covariant derivatives are with respect to $g(t)$, $t \in [t_1, t_2]$.

Now, for each $k \ge 0$, let $\lambda_k=\lambda_k(l)$ be a smooth, compactly supported function on
$[0, \infty)$ such that
$\lambda_k(l)=1, \quad l \in [0, k]$; $0 \le \lambda_k(l) \le 1, \quad l \in [k, k+1]$; and
$\lambda_k(l)=0, \quad l \in [k+1, \infty)$. We also require $|\lambda'_k| \le 4$.
Finally, we take $\phi_k= \lambda_k(L(x))$ as a test function.

 By (\ref{H*Pu}), we have, since $\phi_k$ is compactly supported,
\[
\al
&\frac{d}{dt} \int_{\M} P(u) \phi_k(x) dg(t)
= \int_{\M} [\partial_t P(u)- R P(u) ]\phi_k(x) dg(t)\\
&= \int_{\M} [\partial_t P(u)- R P(u) +\Delta P(u) ]\phi_k(x) dg(t) -
 \int_{\M} P(u)  \Delta \phi_k(x) dg(t)\\
 &=\int_{\M} 2 \tau \left| Ric - Hess \ln u -\frac{g}{2 \tau} \right|^2 u \phi_k(x) dg(t)
 - \int_{\M} P(u)   \Delta \phi_k(x) dg(t).
 \eal
\] Let $t_3, t_4 \in [t_1, t_2]$. Integration on the above yields
\[
\al
&\int^{t_4}_{t_3} \int_{\M} 2 \tau \left| Ric - Hess \ln u -\frac{g}{2 \tau} \right|^2 u \phi_k(x) dg(t)dt\\
&= \int_{\M} P(u) \phi_k(x) dg(t_4) - \int_{\M} P(u) \phi_k(x) dg(t_3)
+ \int^{t_4}_{t_3} \int_{\M}  P(u)  \Delta \phi_k(x) dg(t) dt.
\eal
\] By (\ref{puiu}), this becomes
\[
\al
\int^{t_4}_{t_3} \int_{\M} &2 \tau \left| Ric - Hess \ln u -\frac{g}{2 \tau} \right|^2 u \phi_k(x) dg(t)dt\\
&= \int_{\M} i(u) \phi_k(x) dg(t_4) - \int_{\M} i(u) \phi_k(x) dg(t_3) \\
&\qquad - 2\tau  \int_{\M}  \Delta u \phi_k(x) dg(t_4) + 2 \tau  \int_{\M} \Delta u \phi_k(x) dg(t_3) \\
&\qquad + \int^{t_4}_{t_3} \int_{\M}  i(u)  \Delta \phi_k(x) dg(t) dt
- 2  \int^{t_4}_{t_3} \int_{\M} \tau  \Delta u  \Delta \phi_k(x) dg(t) dt.
\eal
\] After integration by parts, we arrive at
\[
\al
&\int^{t_4}_{t_3} \int_{\M} 2 \tau \left| Ric - Hess \ln u -\frac{g}{2 \tau} \right|^2 u \phi_k(x) dg(t)dt\\
&= \int_{\M} i(u) \phi_k(x) dg(t_4) - \int_{\M} i(u) \phi_k(x) dg(t_3) + \int^{t_4}_{t_3} \int_{\M}  i(u)  \Delta \phi_k(x) dg(t) dt \\
&\qquad - 2\tau  \int_{\M}   u \Delta \phi_k(x) dg(t_4) + 2 \tau  \int_{\M}  u \Delta \phi_k(x) dg(t_3)
- 2  \int^{t_4}_{t_3} \int_{\M} \tau   u  \Delta \Delta \phi_k(x) dg(t) dt.
\eal
\] Notice that the support of $\Delta \phi_k$ is in the region $\{ k \le L(x) \le k+1 \}$.
Since $L(x)$ is comparable with the distance function $d(x, 0, g(t))$, the classical volume
comparison theorem tells us that
$| \{ k \le L(x) \le k+1 \}|_{g(t)} \le C e^{c k}$. Now, recall from (\ref{ubounds}) and
(\ref{iu<}) that $u$ and $i(u)$ have  quadratic exponential decay property.  Also
(\ref{LLLL}) implies that $|\Delta \phi_k| \le C$ and $|\Delta \Delta \phi_k| \le C$. So we can take
$\lim_{k \to \infty}$ inside the integrals on the right hand side the last identity.
On the other hand, $\phi_k$ is a nondecreasing function of $k$, which converges to $1$ pointwise.  Therefore we can
apply the monotone convergence theorem on the left hand side. Therefore
\[
\int^{t_4}_{t_3} \int_{\M} 2 \tau \left| Ric - Hess \ln u -\frac{g}{2 \tau} \right|^2 u  dg(t)dt\\
= \int_{\M} i(u)  dg(t_4) - \int_{\M} i(u)  dg(t_3) .
\] This proves the Corollary.
\qed
\\

Finally, we partially extend Perelman's shrinking breather theorem to the noncompact case.

\begin{definition}  (Breathers)
\lab{debreather}

A Ricci flow $(\M, g(t))$ is a called a breather if for some $t_1<t_2$ and $c >0$ there is the relation $c \phi^* g(t_1)= g(t_2)$
for a diffeomorphism $\phi$. The flow in cases $c=1$, $c<1$, $c>1$ are called steady, shrinking and expanding breathers
respectively.
\end{definition}

When $\M$ is compact, Perelman \cite{P:1} proved that a breather is a gradient Ricci soliton, i.e. the Ricci curvature is
given by the Hessian of a scalar function. For the noncompact case, we have

\begin{proposition}
\lab{prnobreather} Let $(\M, g(t))$ be a noncompact Ricci flow
with bounded geometry in the time interval $[0, T]$. Suppose $(\M,
g(t))$ is a shrinking breather in the sense that  $c \, \phi^*
g(t_1)= g(t_2)$ for some  diffeomorphism $\phi$, $c<1$ and
$t_1<t_2$ where $t_1, t_2 \in (0, T)$.  Suppose also $\mu(g(t_2),
\frac{c (t_2-t_1)}{1-c})<\mu_\infty(g(t_2), \frac{c
(t_2-t_1)}{1-c})$. Then $(\M, g(t))$ is a gradient shrinking
soliton on the time interval $[t_1, T]$. \proof
\end{proposition}

We follow the same strategy as Perelman's proof for the compact case. The new input is the existence of extremal for the $\mu$ invariant in the noncompact setting. Since $\M$ has bounded geometry, we have shown in the proof of Theorem \ref{thmain} (a) that the
Log Sobolev functional is bounded from below by a negative constant. By (\ref{W=L}), the $W$ entropy
also has a lower bound for any finite parameter $\tau$. Thus $\mu(g, \tau)$ is a finite number.

Define $L = \frac{t_2 - c t_1}{1-c}$ where $c$ is the number given in the statement of the proposition.
Then $c (L-t_1)= L-t_2$.  By the scaling and diffeomorphism invariance of the $\mu$ invariant, we have
\[
\mu(g(t_2), L-t_2)=\mu(g(t_2), c(L-t_1))=\mu(c g(t_1), c(L-t_1))=\mu( g(t_1), L-t_1).
\]Note that $L-t_2 = \frac{c (t_2-t_1)}{1-c}$. By the condition $\mu(g(t_2), \frac{c (t_2-t_1)}{1-c})<\mu_\infty(g(t_2), \frac{c (t_2-t_1)}{1-c})$, we can apply Theorem  \ref{thmain} to conclude that $\mu(g(t_2), L-t_2)$ is reached by an extremal
function $v_2$.

Let $u$ be the solution of the
final value problem of the conjugate heat equation:
\[
\begin{cases}
\Delta u - R u + u_t=0, \quad t \in [t_1, t_2]\\
u(x, t_2) = v^2_2\\
\partial_t g(t) = - 2 Ric, \quad t \in [t_1, t_2].
\end{cases}
\]
Since $[t_1, t_2] \subset (0, T)$, Shi's derivative estimate
\cite{Sh:1} shows that the 4-th order derivatives of the curvature
tensor are uniformly bounded in ${\M} \times [t_1, t_2]$. This
allows us to use the Corollary. Let $v=v(x, t) = \sqrt{u(x, t)}$. Since $v_2$ is an extremal of the
W entropy at $t_2$, we know from the Corollary that
\[
\al
&\mu(g(t_2), L-t_2) = W(g(t_2), v(\cdot, t_2), L-t_2) \\
&=  W(g(t_1), v(\cdot, t_1), L-t_1) + \int^{t_2}_{t_1} \int_{\M} 2
\tau  \left| Ric -Hess \ln u - \frac{1}{2 \tau} g \right|^2 \ u \
dg(t) dt. \eal
\] Using $W(g(t_1), v(\cdot, t_1), L-t_1) \ge \mu( g(t_1), L-t_1) =\mu(g(t_2), L-t_2)$, we see that
\[
 \int^{t_2}_{t_1} \tau \int_{\bf M} \left| Ric -Hess \ln u - \frac{1}{2 \tau} g \right|^2
\ u \  dg(t) dt \le 0
\] which implies that $Ric -Hess \ln u - \frac{1}{2 \tau} g =0$. i.e. the Ricci flow is a
 gradient shrinking soliton in the time interval $[t_1, t_2]$. By
 the uniqueness theorem of Chen and Zhu \cite{CZ:1} in the noncompact case,
 the Ricci flow is a gradient shrinking soliton on $[t_1, T]$. This proves the Proposition.
\qed

{\bf Acknowledgement.} I wish to thank Professor Zhiqin Lu for a useful conversation and
Professor Bennett  Chow and Professor Lei Ni for their continuous help on my studying of Ricci flow over the years.

\bigskip

\noindent e-mail:  qizhang@math.ucr.edu

\enddocument